\documentclass{amsart}
\usepackage{latexsym,amssymb,amscd}

\theoremstyle{plain} 

\newtheorem{theorem}{Theorem}[section]

\newtheorem{lemma}{Lemma}[section]

\newtheorem{proposition}{Proposition}[section]

\theoremstyle{remark} 
\newtheorem{remark}{Remark}[section] 
\theoremstyle{remark}

\theoremstyle{definition}

\newtheorem{definition}{Definition}[section]
\newtheorem{example}{Example}[section]
\numberwithin{equation}{section}
\def\<{\left < }
\def\>{\right >} 
\def\({\left ( } 
\def\){\right )} 
\begin{document}

\title[Segre embedding  in
differential geometry] {Segre embedding and related maps and
immersions in differential geometry}

\author{Bang-Yen Chen}

 \address{Department of Mathematics\\
	Michigan State University\\ East Lansing, MI 48824--1027, USA}

\email{bychen@math.msu.edu}

\begin{abstract}  Segre embedding was
introduced  by C. Segre (1863--1924) in his famous 1891
article \cite{segre}. The Segre embedding plays an important roles
in  algebraic geometry as well as in differential geometry,
mathematical physics, and coding theory. In this article, we  survey
main results on Segre embedding in differential geometry. Moreover,
we also present recent differential geometric results on maps and
immersions which are constructed in ways similar to Segre embedding. 
\end{abstract} 



\insert\footins{\noindent \footnotesize
  Arab. J. Math. Sci. 8 (2002), 1--39.}

\maketitle
\baselineskip.18in

\centerline{{\sc Table of Contents.}\footnote{The author would
like to express his thanks to the editors for the invitation to
publish this article in this Journal}} 
\vskip.2in
\begin{enumerate}
 \item[1.]  Introduction.

\item[2.]  Basic formulas and definitions.

\item[3.] Differential geometric characterizations of
Segre embedding.

\item[4.] Degree of K\"ahlerian immersions and homogeneous
K\"ahlerian submanifolds via Segre embedding.

\item[5.] $CR$-products and Segre embedding.

\item[6.] $CR$-warped products and partial
Segre $CR$-immersions.

\item[7.] Real hypersurfaces as partial Segre
embeddings.

\item[8.] Complex extensors, Lagrangian submanifolds and Segre
embedding.

\item[9.] Partial Segre $CR$-immersions in complex projective
space. 

\item[10.] Convolution of Riemannian manifolds.

\item[11.] Convolutions and Euclidean Segre maps.

\item[12.] Skew Segre embedding. 

\item[13.] Tensor product immersions and Segre embedding.

\item[14.] Conclusion.

\item[]  References
\end{enumerate}

\section{Introduction.} 

Throughout this article, we denote by $\,CP^n(c)\,$ the complex
projective $n$-space endowed with the Fubini-Study metric with
constant holomorphic sectional curvature c. Let
$(z_0,\ldots,z_n)$ denote a homogeneous coordinate system on
$CP^n(c)$

There are two well-known examples of algebraic manifolds in
complex projective spaces: the Veronese embedding
$v_n: CP^n(2)\to CP^{n(n+3)/2}(4)$ and the Segre
embedding $S_{hp}:CP^{h}(4)\times CP^{p}(4)\to CP^{h+p+hp}(4)$. 

The Veronese embedding $v_n$ is a K\"ahlerian embedding, that
is a holomorphically isometric embedding, of $CP^n(2)$ into
$CP^{n(n+3)/2}(4)$ given
by  homogeneous monomials of degree 2:
\begin{align} \label{1.0} &
v_{n}:CP^{n}(2)\to CP^{n(n+3)/2}(4);\,
\\& (z_0,\ldots,z_n)\mapsto
\Bigg(z_0^2,\sqrt{2}z_0
z_1,\ldots,\sqrt{\frac{2}{\alpha_i !\alpha_j
!}}z_i^{\alpha_i}z_j^{\alpha_j},\ldots, z_n^2\Bigg)\notag\end{align}
with $\alpha_i+\alpha_j=2$.  For $n=1$, this is nothing but the
quadric curve $$Q_1=\Bigg\{(z_0,z_1,z_2)\in CP^2: \sum_{j=0}^2
z_j^2=0\Bigg\}$$ in $CP^2(4)$.
 
The Veronese embedding can be extended to $\alpha$-th Veronese
embedding $v_n^\alpha$ with $\alpha\geq 2$: 
\begin{align}v_n^\alpha:CP^n\Bigg(\frac{4}{\alpha}\Bigg)\to
CP^{\binom{n+\alpha}  {\alpha}-1}(4)
\end{align}
defined by
\begin{align}  (z_0,\ldots,z_n)\mapsto
\Bigg(z_0^{\alpha},\sqrt{\alpha}z_0^{\alpha-1}
z_1,\ldots,\sqrt{\frac{\alpha}{\alpha_0 !\cdots\alpha_n
!}}z_0^{\alpha_0}\cdots z_n^{\alpha_n},\ldots, z_n^{\alpha}\Bigg)
\end{align}
with $\alpha_0+\cdots+\alpha_n=\alpha$.

On the other hand, the Segre embedding:
\begin{align} \label{1.1}
S_{hp}:CP^{h}(4)\times CP^{p}(4)\to CP^{h+p+hp}(4),\end{align}
is defined by \begin{align} \label{1.2}
S_{hp}(z_0,\ldots,z_{h};w_0,\ldots,w_p)=\big(z_jw_t \big)_{0\leq
j\leq h,0\leq t\leq p},\end{align} 
where  $(z_0,\ldots,z_{h})$ and
$(w_0,\ldots,w_p)$ are the homogeneous coordinates of
$CP^{h}(4)$ and $CP^p(4)$, respectively. 
This embedding \eqref{1.1} was introduced
by C.  Segre in 1891 (see \cite{segre}). It is well-known 
that the Segre embedding $S_{hp}$ is also a K\"ahlerian
embedding. 

When $h=p=1$,
the Segre embedding is nothing but the complex quadric
surface,
$Q_2=CP^1\times CP^1$ in $CP^3$, defined by
\begin{align} 
Q_2=\Bigg\{(z_0,z_1,z_2,z_3)\in CP^3: \sum_{j=0}^3
z_j^2=0\Bigg\}.\end{align}

The Segre embedding can also be naturally extended to product
embeddings of arbitrary number of complex projective spaces
as follows.

Let $(z_0^i,\ldots,z_{n_i}^i)$ $(1\leq i\leq
s)$ denote the homogeneous coordinates of
$CP^{n_i}$. Define a map:
\begin{align} \label{1.4}
S_{n_1\cdots
n_s}:CP^{n_1}(4)\times\cdots\times CP^{n_s}(4)\to
CP^N(4),\end{align} $$ N=\prod_{i=1}^s (n_i+1)-1,$$ which
maps each point
$((z_0^1,\ldots,z_{n_1}^1),\ldots,
(z_0^s,\ldots,z_{n_s}^s))$ in the product
K\"ahlerian manifold
$CP^{n_1}(4)\times\cdots\times CP^{n_s}(4)$ to
the point $(z^1_{i_1}\cdots
z^s_{i_j})_{1\leq i_1\leq
n_1,\ldots,1\leq i_s\leq n_s}$ in
$CP^N(4)$. This map
$S_{n_1\cdots n_s}$ is also a K\"ahlerian embedding.

The Segre embedding is known to be the
simplest K\"ahlerian embedding from product algebraic
manifolds into complex projective spaces. It is well-known
that the Segre embedding plays an important role in algebraic
geometry (see \cite{mumford,Sm}). The Segre embedding has
also been applied to  differential geometry as well as to
coding theory (see, for instance,
\cite{In,Sc,Sk}) and to mathematical physics (see, for
instance,
\cite{br,Se}).  

The purpose of this article is to survey the main results on
Segre embedding in differential geometry. Furthermore, we
also present recent results in  differential
geometry concerning maps, immersions, and embedding which are
constructed in ways similar to Segre embedding defined by
\eqref{1.2}.

\section{Basic formulas and definitions.}

Let $M$ be a Riemannian
$n$-manifold with inner product
$\<\;\, ,\;\>$ and let
$e_1,\ldots,e_n$ be an
orthonormal frame fields on $M$.
For a differentiable function
$\varphi$ on $M$, the gradient
$\nabla\varphi$ and the Laplacian
$\Delta\varphi$ of $\varphi$ are
defined respectively by
\begin{align}\label{E:2.1} &\<\nabla
\varphi,X\>=X\varphi,
\\ \label{E:2.2}& \Delta\varphi=\sum_{j=1}^n
\big\{e_je_j\varphi -(\nabla_{e_j}e_j)\varphi
\big\}\end{align}for
vector fields $X$ tangent to $M$,
where $\nabla$ is the
Levi-Civita connection on $M$. 

If $M$ is
isometrically immersed in a
Riemannian manifold $\tilde M$, then the formulas of
Gauss and Weingarten for $M$
in $\tilde M$ are given respectively by
\begin{align}\label{E:2.3} &\tilde \nabla_XY=\nabla_X
Y+\sigma(X,Y),\\ \label{E:2.4} &\tilde
\nabla_X\xi =-A_\xi X+D_X\xi\end{align} for vector
fields $X,\,Y$ tangent to
  $N$ and  $\xi$ normal to $M$,
where $\tilde \nabla$ denotes
the Levi-Civita connection
on $\tilde M$, $\sigma$ the second fundamental form, $D$
the normal connection, and $A$
the  shape operator of $M$ in $\tilde M$. 

 The second
fundamental form and the shape
operator are related by
$$\<A_\xi X,Y\>=\<\sigma(X,Y),
\xi\>.$$

The mean curvature vector $\overrightarrow H$ is defined by
\begin{align} \overrightarrow H = {1\over n}\,\hbox{\rm
trace}\,\sigma = {1\over n}\sum_{i=1}^{n}
h(e_{i},e_{i}), \end{align}
where $\{e_{1},\ldots,e_{n}\}$ is a local  orthonormal
frame of the tangent bundle $TN$ of $N$. The
{\it squared mean curvature\/} is given by
$H^2=\left<\right.\overrightarrow H,\overrightarrow
H\left.\right>$, where $\<\;\,,\;\>$ denotes the inner
product. A submanifold
$N$  is called {\it totally geodesic\/}  in $\tilde
M$ if the second fundamental form  of $N$ in $\tilde M$
vanishes identically. And $N$ is called {\it minimal} if its mean
curvature vector vanishes identically.

The  {\it equation of
Gauss\/} is given by
\begin{align}\label{E:2.5} {\tilde
R}(X,Y;Z,W)=&\,R(X,Y;Z,W)+\<\sigma(X,Z),\sigma(Y,W)\>\\
&-\<\sigma(X,W),\sigma(Y,Z)\>,\notag \end{align} for
$X,Y,Z,W$ tangent to $M$, where $R$ and $\tilde R$
denote the curvature tensors of $M$ and
$\tilde M$, respectively.

For the second fundamental form $\sigma$, we
define its covariant derivative ${\bar
\nabla}\sigma$ with respect to the connection
on
$TM \oplus T^{\perp}M$ by
\begin{equation}\label{E:2.6}({\bar\nabla}_{X}\sigma)(Y,Z)=D_{X}(\sigma
(Y,Z))-\sigma(\nabla_{X}Y,Z)
-\sigma(Y,\nabla_{X}Z)\end{equation} for $X,Y,Z$ tangent to
$M$. The {\it equation of Codazzi\/} is
\begin{equation}\label{E:2.7}({\tilde
R}(X,Y)Z)^{\perp}=({\bar\nabla}_{X}\sigma)(Y,Z)-
({\bar\nabla}_{Y}\sigma)(X,Z),\end{equation} where $({\tilde
R}(X,Y)Z)^{\perp}$ denotes the
normal component of ${\tilde
R}(X,Y)Z$.  

A submanifold $M$ in a Riemannian manifold $\tilde M$ is said
to have {\it parallel second fundamental form} if $\bar\nabla
\sigma=0$ identically.

 The Riemann curvature tensor of a complex space form
$\tilde M^m(4c)$ of constant holomorphic sectional
curvature $4c$ is given by
\begin{align}\label{E:2.9} \tilde R(X,Y;Z,W)=\; &
c\, \big\{\<X,W\>\<Y,Z\>-\<X,Z\>\<Y,W\>+\<JX,W\>\<JY,Z\>
\\&-\<JX,Z\>\<JY,W\>+2\<X,JY\>\<JZ,W\>
\big\}.\notag\end{align} 

If we define the $k$-th ($k\geq 1$) covariant derivative of the
second fundamental form $\sigma$ by
\begin{align}\label{E:2.10} \big(\bar\nabla^k
\sigma\big)(X_1,&\ldots,X_{k+2})=D_{X_{k+2}}\big(\big(\bar
\nabla^{k-1}\sigma\big)(X_1,\ldots,X_{k+1})\big)\\&
-\sum_{i=1}^{k+1} (\big(\bar
\nabla^{k-1}\sigma\big)(X_1,\ldots,\nabla_{X_{k+2}}X_i,\ldots,X_{k+1}),
\notag\end{align}
then $\bar\nabla^k \sigma$ is a normal-bundle-valued tensor of type
$(0,k+2)$. Moreover, it can be proved that $\bar\nabla^k \sigma$
satisfies
\begin{align}\label{E:2.11} \big(\bar\nabla^k
\sigma\big)(X_1,&X_2,X_3,\ldots,X_{k+2})-\big(\bar\nabla^k
\sigma\big)(X_2,X_1,X_3,\ldots,X_{k+2})\\&=R^D(X_1,X_2)
\big(\big(\bar
\nabla^{k-2}\sigma\big)(X_3,\ldots,X_{k+2})\big)\notag\\&
+\sum_{i=3}^{k+2} (\big(\bar
\nabla^{k-2}\sigma\big)(X_3,\ldots,R(X_1,X_2)X_i,\ldots,X_{k+2}),
\notag\end{align}
for $k\geq 2$. For simplicity, we put
$\bar\nabla^0\sigma=\sigma$.

\section{Differential geometric characterizations of Segre
embedding.}

H. Nakagawa and R. Tagaki \cite{na} classify
complete K\"ahlerian submanifold in complex projective spaces
with parallel second fundamental form. In particular, they show
that the Segre embedding are the only one which are reducible.
More precisely, they prove the following.

\begin{theorem} Let $M$ be a complete K\"ahlerian submanifold of $n$
complex dimensions embedded in $CP^m(4)$.  If $M$ is reducible
and has parallel second fundamental form,
then $M$ is congruent to
$CP^{n_1}(4)\times CP^{n_2}(4)$ with $n=n_1+n_2$ and the
embedding is given by the Segre embedding.\end{theorem}

The following results published in 1981 by B. Y. Chen 
\cite{c1} for $s=2$ and
by B. Y. Chen and W. E. Kuan \cite{ck1,ck2} for $s\geq 3$ can
be regarded as ``converse'' to Segre embedding constructed in
1891 by C. Segre.  

\begin{theorem}[\cite{c1,ck1,ck2}]\label{T:s1} Let
$M_1^{n_1},\ldots,M_s^{n_s}$ be K\"ahlerian manifolds of
dimensions
$n_1,\ldots,$ $n_s$, respectively. Then locally every
K\"ahlerian immersion
$$f:M_1^{n_1}\times\cdots\times M_s^{n_s}\to
CP^N(4),\quad N=\prod_{i=1}^s (n_i+1)-1,$$ of
$M_1^{n_1}\times\cdots\times M_s^{n_s}$ into $CP^N(4)$  
is   the Segre embedding, that is,
$M_1^{n_1},\ldots,M_s^{n_s}$ are open portions of
$CP^{n_1}(4),\ldots, CP^{n_s}(4)$, respectively.
Moreover, the K\"ahlerian immersion $f$ is given by
the Segre embedding.\end{theorem}

Let $||\bar\nabla^k\sigma||^2$ denote  the squared norm of the
$k$-th covariant derivative of the second fundamental form. The
Segre embedding can also be characterized by 
$||\bar\nabla^k\sigma||^2$ as given in the following.

 \begin{theorem}[\cite{c1,ck1,ck2}]\label{T:s2} Let
$M_1^{n_1}\times\cdots\times M_s^{n_s}$ be a
product K\"ahlerian submanifold in 
$CP^m(4)$ of arbitrary codimension. Then we have
\begin{equation}\label{E:6.2} ||\bar
\nabla^{k-2}\sigma||^2\geq k!\,2^k\sum_{i_1<\cdots<i_k}
n_1\cdots n_k,\end{equation} for $k=2,3,\cdots.$

The equality sign of \eqref{E:6.2}  holds for some $k\geq 2$
if and only if  
$M_1^{n_1},\ldots,M_s^{n_s}$ are open portions of 
$CP^{n_1}(4),\ldots, CP^{n_s}(4)$, respectively, and
 the K\"ahlerian immersion is given by
the Segre embedding.\end{theorem}

When $k=2$, Theorem \ref{T:s1} and Theorem \ref{T:s2}
reduce to

\begin{theorem}[\cite{c1}]\label{T:s3} Let
$M_1^{h}$ and $M_2^{p}$ be two K\"ahlerian manifolds of complex
dimensions
$h$ and $p$, respectively. Then every K\"ahlerian
immersion
$$f:M_1^{h}\times M_2^{p}\to
CP^{h+p+hp}(4),$$ of
$M_1^{h}\times M_2^{p}$ into $CP^{h+p+hp}(4)$  
is locally  the Segre embedding, that is,
$M_1^{h}$ and $M_2^p$ are open portions of
$CP^{h}$ and $ CP^{p}$, respectively, and
moreover, the K\"ahlerian immersion $f$ is given by
the Segre embedding.\end{theorem}
 
\begin{theorem}[\cite{c1}]\label{T:s4} Let
$M_1^{h}\times M_2^{p}$ be a product K\"ahlerian
submanifold  in 
$CP^m(4)$ of arbitrary codimension. Then we have
\begin{equation}\label{E:3.2} ||\sigma||^2\geq
8hp.\end{equation}

The equality sign of \eqref{E:3.2}  holds if and
only if  $M_1^h$ and $M_2^p$ are open portions of
$CP^{h}(4)$ and $CP^{p}(4)$, respectively, and
 the K\"ahlerian immersion is given by
the Segre embedding $S_{h,p}$.
\end{theorem}

Let $CP^n_s(4)$ denote the indefinite complex projective
space of complex dimension $n$, index $2s$, and constant
holomorphic sectional curvature $4$ and let
$S^{2n+1}_{2s}(1)$  be the $(2n+1)$-dimensional indefinite
unit-sphere with index $2s$ and of constant sectional
curvature 1.
Thus a point of 
$CP^n_s(c)$ can be represented by $[(z,w)]$, where
$z=(z_1,\ldots,z_s)\in {\bf C}^s,\,
w=(w_1,\ldots,w_{n-s-1})\in {\bf C}^{n-s+1}$, $(z,w)\in
S^{2n+1}_{2s}(1)\subset {\bf C}^{n+1}_s$ and  $[(z,w)]$ is
the equivalent class of the Hopf projection
$$\pi_H:S^{2n+1}_{2s}(1)\to CP^n_s(4).$$

Consider the map:
$$\phi:CP^h_s(4) \times CP^p_t(4)\to 
CP^{N(h,p)}_{R(h,p,s,t)}(4)$$
with 
\begin{align} &N(h,p)=h+p+hp,\notag\\
 R(h,p,s,& t)=s(p-t)+t(h-s)+s+t\notag
\end{align}
given by
$$ \phi([(z,w)],[(x,y)])=[(x_i y_\alpha, w_k x_a,z_j
x_b,w\ell y_\beta)]$$
for $$1\leq i,j\leq s; \, 1\leq k,\ell\leq h-s+1;\, 1\leq
a,b\leq t;\, 1\leq \alpha,\beta \leq h-t+1.$$ 

Then map $\phi$ is
a well-defined holomorphic isometric embedding, which is
called the indefinite Segre embedding $CP^h_s(4) 
\times CP^p_t(4)$ into $CP^{N(h,p)}_{R(h,p,s,t)}(4)$. 

T. Ikawa, H. Nakagawa and A. Romero  \cite{ikawa} study
indefinite version of Theorem \ref{T:s3} and obtain the
following.

\begin{theorem} Let $M_s^n$ and
$M_t^m$ be two complete indefinite K\"ahlerian manifolds with
complex dimensions
$n$ and
$m$, and indices $2s$ and $2t$, respectively. If there exists a
holomorphic isometric immersion from the product $M_s^n\times
M_t^m$ into an indefinite complex projective space of
complex dimension $N$ and index $2r$, then we have:

\begin{enumerate}
\item $N\geq
n+m+nm$ and $ r\geq s(m-t)+t(n-s)+s+t$.

\item If $N=n+m+nm$, then
the immersion is obtained by the indefinite Segre imbedding.
\end{enumerate}\end{theorem}

\begin{remark} The assumption of ``K\"ahlerian immersion'' in
Theorems \ref{T:s1}, \ref{T:s2}, and \ref{T:s3} is necessary. In
fact, let
$M_i$ be a projective nonsingular embedded variety of dimension
$n_i\geq 1\, (i=1,2,\cdots,r)$ and let
$$M=M_1\times \cdots\times M_r\subset
CP^{n_1}\times\cdots\times CP^{n_r}
\xrightarrow[\text{{\rm
Segre embedding}}]{\text{$S_{n_1\cdots n_r}$}}
CP^N$$ be the composition embedding from the product
$M_1\times \cdots\times M_r$ into $CP^N$ with
$N=\prod_{i=1}^r (n_i+1)-1$ via the Segre embedding. 

M. Dale
considers in \cite{D} the problem of finding the embedding
dimension $e$ such that $M$ can be  embedded (not necessary
K\"ahlerian embedded in general) in
$CP^e$, but not in $CP^{e-1}$, via a projection. Using an
algebraic result of A. Holme, Dale characterizes $e$ in
terms of the degree of the Segre classes of $M$, he proves that
$e=2(n_1+\cdots+n_r)+1$, unless $r=2$, $X_1=CP^{n_
1}$, $X_2=CP^{n_2}$, in which case $e=2(n_1+n_2)-1$.
\end{remark}

\section{Degree of K\"ahlerian immersions and homogeneous
K\"ahlerian submanifolds via Segre embedding.}

 By {\it applying Segre embedding}, R. Takagi and M. Takeuchi
define in  \cite{takagi} the notion of tensor products of
K\"ahlerian immersions in complex projective spaces as
follows:

Suppose that $f_i:M_i\to CP^{N_i}(4),\; 
i=1,\ldots,s,$ are full K\"ahlerian embeddings
of irreducible Hermitian symmetric spaces of
compact type.  Consider the composition given by
\begin{align}f_1\boxtimes\cdots\boxtimes
f_s:M_1\times\cdots&\times  M_s\xrightarrow[\text{{\rm
product embedding}}]{\text{$f_1\times\cdots\times
f_s$}} CP^{N_1}\times \cdots\times CP^{N_s}
\\& \xrightarrow[\text{{\rm
Segre embedding}}]{\text{$S_{N_1\cdots
N_s}$}}CP^N(4),\notag\end{align}  with $N=\prod_{i=1}^s
(N_i+1)-1$. This composition is a full K\"ahlerian
embedding, which is called the {\it tensor product} of
$f_1,\ldots,f_s$.
 
 H. Nakagawa and R. Tagaki in
\cite{na} and R. Tagaki and M. Takeuchi in \cite{takagi} 
had obtained a close relation between the
{\it degree} and the {\it rank} of a symmetric K\"ahlerian
submanifold in complex projective space;
namely, they proved the following.

\begin{theorem} Let $f_i:M_i\to CP^{N_i}(4),\; i=1,\ldots,s,$
are
$p_i$-th full K\"ahlerian embeddings of irreducible
Hermitian symmetric spaces of compact type. Then
the degree of the tensor product $f_1\boxtimes\cdots\boxtimes
f_s$ of
$f_1,\ldots,f_s$ is given by $\sum_{i=1}^s r_ip_i$,
where $r_i=rank(M_i)$. \end{theorem}

Related with this theorem, we mention the following result by M.
Takeuchi \cite{T} for K\"ahlerian immersions of homogeneous
K\"ahlerian manifolds.

\begin{theorem} Let $f:M\to CP^m(4)$ be a K\"ahlerian
immersion of a globally homogeneous K\"ahlerian manifold $M$.
Then

\begin{enumerate}
\item $M$ is compact and simply-connected;

\item $f$ is an embedding; and 

\item $M$ is the orbit in $CP^m(4)$ of the highest
weight in an irreducible unitary
representation of a compact semisimple Lie
group. 

\end{enumerate}\end{theorem}

The notion of the {\it degree} of K\"ahlerian
immersions in the sense of 
\cite{takagi} is defined as follows:
Let $V$ be a real vector space of dimension $2n$
with an almost complex structure $J$ and a Hermitian inner
product $g$. Denote the complex linear extensions of $J$ and
$g$ to the complexification $V^C$ of $V$ by the same $J$ and
$g$, respectively. Let $V^+$ and $V^-$ be the eigensubspace
of $J$ on $V^C$ with eigenvalue $1$ and $-1$, respectively.
Then
$V^C=V^+\oplus V^-$ is an orthogonal direct sum with respect
to the inner product. 

Let $E$ be a real vector bundle over a manifold $M$ with
a Hermitian structure $(J,g)$ on fibres.  
The Hermitian structure induces a  Hermitian inner
product on  $E^C$. 
We  have subbundle $E^+$ and $E^-$ such that $E^C$
satisfies
$E^C=E^+\oplus E^-$ and the
complex conjugation $E^\pm \xrightarrow{\text{{\rm --}}}
E^\mp$. The map on the space of sections induced from
the complex conjugation is denoted by
$\Gamma (E^\pm)\xrightarrow{\text{{\rm --}}}\Gamma (E^\mp).$
Let $(M,g,J)$ be a K\"ahlerian manifold. 
Then we get a Hermitian inner product on the
complexification $T(M)^C$ and subbundles $T(M)^\pm$ of
$T(M)^C$ such that 
$T(M)^C=T(M)^+\oplus T(M)^-$. 
 
 Let  $f:(M,g,J)\to (M',g',J')$
be a K\"ahlerian immersion between K\"ahlerian manifolds.
 The Levi-Civita
connections of $M$ and $M'$ are denoted by $\nabla$ and
$\nabla'$. The induced bundle $f^*T(M')$ has a Hermitian
structure $(J',g')$ induced from the one on $M'$. Also it has
a connection $\nabla'$ induced from  $M'$. 

If we denote the orthogonal complement of $f_*T(x(M)$ in
$T_{f(x)}(M')$ by $T^\perp_x(M)$, then
 $T^\perp(M)$ is a subbundle of
$f^*T(M')$.  
We have the orthogonal Whitney sum decompositions:
$f^*T(M')=f_*T(M)\oplus T^\perp(M),\,
f^*T(M')^C=f_*T(M)^C\oplus T^\perp(M)^C,
$ and $f^*T(M')^\pm=f_*T(M)^\pm\oplus T^\perp(M)^\pm$.

 The orthogonal projection
$f^*T(M')\to T^\perp(M)$ is denoted by $X\mapsto
X^\perp$ and the induced projection
$\Gamma(f^*T(M'))\to \Gamma(T^\perp(M))$ is denoted
by $\xi\mapsto \xi^\perp$.
The normal connection $D$ on $T^\perp(M)$ satisfies
$D_X\xi=(\nabla^{'}_X \xi)^\perp$.

Let $\sigma$ denote the second fundamental
form of $f$. We have
\begin{align}\label{E:4.2.1}\sigma(T_x(M)^+,T_x(M)^-)=\{0\},
\quad \sigma(T_x(M)^\pm,T_x(M)^\pm)\subset
T^\perp_x(M)^\pm.\end{align}

Put $\sigma_2=\sigma$. For $k\geq 3$, we define $\sigma_k$
inductively just like
\eqref{E:2.10} by
\begin{align}\label{E:4.2.2} \sigma_{k+1}&
(X_1,\ldots,X_{k+1})=D_{X_{k+1}}
\sigma_{k}(X_1,\ldots,X_{k})\\&
-\sum_{i=2}^{k}
\sigma_k(X_1,\ldots,\nabla_{X_{k+1}}X_i,\ldots,X_{k}),
\notag\end{align}
for $X_i\in T_x(M)$.

Equations \eqref{E:4.2.1} and \eqref{E:4.2.2} imply that
$\sigma_k(X_1,\ldots,X_k)\in T^\perp_x(M)^+$ for
$X_1,X_2\in T_x(M)^+$ and $X_3,\ldots,X_k\in T_x(M)^C$.

Let $H^k\in \Gamma(\hbox{\rm
Hom} (\otimes^k T(M)^+,T^\perp(M)^+))$ $(k\geq 2)$ be defined
by
$$H^k(X_1,\ldots,X_k)=\sigma_k(X_1,\ldots,X_k),\quad X_i\in
T_x(M)^+.$$
We put $$h=:\sum_{k\geq 2} \sigma_k \in \Gamma\Big( \hbox{\rm
Hom} \Big(\sum_{k\geq 2}\otimes^k
T(M),T^\perp(M)\Big)\Big), $$ 
$$H=:\sum_{k\geq 2} H^k \in \Gamma\Big( \hbox{\rm
Hom} \Big(\sum_{k\geq 2}\otimes^k
T(M)^+,T^\perp(M)^+\Big)\Big).
$$ 
For an integer $k>0$, we define a subspace $\mathcal
H^k_x(M)$ of $T_{f(x)}(M')^+$ to be the subspace spanned by
$T_x(M)^+$ and $H\big(\sum_{2\leq j\leq k} \otimes^j T_x(M)^+ 
\big)$. Then we get a series:
$${\mathcal H}^1_x(M)\subset {\mathcal H}^2_x(M)\subset
\cdots\subset {\mathcal H}^k_x(M)\subset {\mathcal
H}^{k+1}_x(M)\subset \cdots
\subset T_{f(x)}(M')^+$$
of increasing subspaces of $T_{f(x)}(M')^+$.
Let $O^k_x(M)$ be the orthogonal complement of $\mathcal
H^{k-1}_x(M)$ in $\mathcal H^k_x(M)$, where $\mathcal
H^0_x(M)$ is understood to be
$\{0\}$. Then we have an orthogonal direct sum:
$\mathcal H^k_x(M)=O^1_x(M)\oplus O^2_x(M)\oplus\cdots\oplus
O^k_x(M).$

Define $\mathcal R_1=M$. For an integer $k>1$, we define the
set $\mathcal R_k$ of {\it $k$-regular} points of $M$
inductively by
$$\mathcal R_k=\Big\{u\in \mathcal R_{k-1}:\dim_{\hbox{\bf C}}
\mathcal H^k_x(M)=\max_{y\in \mathcal R_{k-1}}\dim_{\hbox{\bf C}}
\mathcal H^k_y(M)\Big\}.$$
Then we have the inclusions:
$\mathcal R_1\supset \mathcal R_2\supset \cdots
\supset\mathcal R_k \supset \mathcal R_{k+1}\supset\cdots.$
Note that each $\mathcal R_k$ is an open nonempty subset 
of $M$ and,  for each $k$,
$\mathcal H^k(M)=\cup_{x\in\mathcal R_k} \mathcal H^k_x(M)$
is a  complex vector bundle over $\mathcal R_k$ which
is a subbundle of $f^*T(M')^+|_{\mathcal R_k}$.

For an integer $k\geq 1$ and a point $x\in \mathcal
R_k$, we have \cite{takagi}

(1) $\nabla'_X Y\in\mathcal H^{k+1}_x(M)$ for  $X\in
T_x(M)^+$ and local sections $Y$ of $\mathcal H^k(M)$;

(2) $O^{k+1}_x(M)=\{0\}$ if and only if, for each $X\in
T_x(M)^+$ and each local  section $Y$ of $\mathcal
H^k(M)$, we have  $\nabla'_X Y\in\mathcal
H^{k}_x(M)$.

Thus, there is a unique integer $d>0$ such that 
$O^d_x(M)\ne \{0\}$ for some $x\in{\mathcal
R}_d$ and $ O^{d+1}_x(M)= \{0\}$ for each
$x\in{\mathcal R}_d$.
The integer $d$ is called the {\it degree} of the K\"ahlerian
immersion $f:M\to \tilde M$.

\section{$CR$-products and Segre embedding.} 

A submanifold $N$ in a K\"ahlerian manifold $\tilde
M$ is called a {\it totally real submanifold} \cite{chenogiue1}
if the complex structure $J$ of $\tilde M$ carries each tangent
space of $N$ into its corresponding normal space, that is,
$JT_x N\subset T_x^\perp N$, $x\in N$.  An $n$-dimensional
totally real submanifold in a K\"ahlerian manifold $\tilde M^n$
with complex dimension $n$ is called a {\it Lagrangian
submanifold}. (For  latest surveys on Lagrangian submanifolds
form differential geometric point of view, see
\cite{c2,2001}).

A submanifold $N$ in a K\"ahlerian manifold $\tilde
M$ is called a {\it $CR$-submanifold} \cite{bejancu} if there
exists on $N$ a holomorphic distribution $\mathcal D$ whose
orthogonal complement $\mathcal D^\perp$ is a totally real
distribution, that is, $J\mathcal D^\perp_x\subset
T^\perp_x N$. 

The notion of $CR$-products was introduced 
 in \cite{c1} as follows: A $CR$-submanifold
$N$ of a K\"ahlerian manifold 
$\tilde M$ is called a {\it $CR$-product} if locally it is
a Riemannian product of a K\"ahlerian
submanifold $N_T$ and a totally real
submanifold $N_\perp$ of $\tilde M$. 

For a $CR$-submanifold
$N$ in a K\"ahlerian manifold $\tilde M$, we put
$$JX=PX+FX,\quad X\in TN,$$
where $PX$ and $FX$ denote
the tangential and the normal components of $JX$, respectively.

It is proved in \cite{c1} that a submanifold $M$ of
a K\"ahlerian manifold is a $CR$-product if and
only if $\nabla P=0$ holds, that is,  $P$ is
parallel with respect to the Levi-Civita
connection of $M$.

\begin{example} Let $\psi_1:N_T\to CP^{n_1}(4)$ be a
K\"ahlerian immersion of a K\"ahlerian manifold $N_T$ into
$CP^{n_1}(4)$ and let
$\psi_2:N_\perp\to CP^{n_2}(4)$ be a totally real immersion of
a Riemannian $p$-manifold $N_\perp$ into $CP^{n_2}(4)$.
Then the composition:
\begin{align} S_{n_1 n_2}\circ (\psi_1,\psi_2):N_T\times N_\perp
&\xrightarrow[\text{\rm product immersion
}] {\text{{$(\psi_1,\psi_2)$
}}} CP^{n_1}(4)\times CP^{n_2}(4)\\&
\xrightarrow[\text{\rm Segre embedding}] {\text{{$S_{n_1
n_2}$ }}} CP^{n_1+n_2+n_1 n_2}(4) \notag
\end{align} is isometric immersed as a $CR$-product in
$CP^{n_1+n_2+n_1 n_2}(4)$.

In particular, if  $\iota:CP^h(4)\to CP^h(4)$ is the
identity map of $CP^h(4)$ and $\varphi:N_\perp\to CP^p(4)$ is a
Lagrangian immersion of a Riemannian $p$-manifold $N_\perp$
into $CP^p(4)$, then the composition:
\begin{align} S_{hp}\circ (\iota,\varphi):CP^h(4)\times N_\perp
&\xrightarrow[\text{\rm product immersion
}] {\text{{$(\iota,\varphi)$
}}} CP^h(4)\times CP^h(4)\\&\xrightarrow[\text{\rm Segre
embedding}] {\text{{$S_{hp}$ }}} CP^{h+p_hp}(4)\notag
\end{align}
is a $CR$-product in $CP^{h+p+hp}(4)$ which is called a {\it standard
$CR$-products} \cite{c1}.
\end{example}

For $CR$-products in complex space forms, the following
results are known. 

\begin{theorem}[\cite{c1}]\label{T:7.1} A
$CR$-submanifold in the complex Euclidean $m$-space $\hbox{\bf
C}^m$  is a $CR$-product if and only if it is a
direct sum of a K\"ahlerian
submanifold and a totally real
submanifold of  linear complex
subspaces. \end{theorem}
 
\begin{theorem} [\cite{c1}] There do not exist
$CR$-products in complex hyperbolic
spaces other than K\"ahlerian submanifolds and
totally real submanifolds.\end{theorem}

 $CR$-products $N_T^h\times N^p_\perp$ in 
$CP^{h+p+hp}(4)
$ are obtained from the
Segre embedding as given in Example 5.1. More precisely, we have the
following.

\begin{theorem}[\cite{c1}]  Let $N_T^{h}\times
N_\perp^p$ a $CR$-product in $CP^m(4)$ with
$\dim_{\hbox{\bf C}} N_T=h$ and
$\dim_{\hbox{\bf R}} N_\perp=p$. Then we have:
\begin{equation}\label{E:7.a} m\geq h+p+ hp
\end{equation}
 
The equality sign of \eqref{E:7.a} holds if and only if the
following statements hold:

$(a)$ $N_T^h$ is an open portion of $CP^h(4)$. 

$(b)$ $N_\perp^p$ is a totally real submanifold. 

$(c)$ The immersion is the following composition:
\begin{equation}\label{E:c} N_T^h\times N_\perp^p
\xrightarrow{\text{{\rm}}}
CP^h(4)\times CP^p(4)\xrightarrow [\text{{\rm
Segre imbedding}}]{\text{$S_{hp}$}} CP^{h+p+hp}(4).
\notag\end{equation} 
\end{theorem}  

\begin{theorem}[\cite{c1}]\label{T:7.3}  Let $N_T^{h}\times
N_\perp^p$ be a $CR$-product in $CP^m(4)$. Then 
 the squared norm of the second
fundamental form  satisfies
\begin{equation}\label{E:7.b} ||\sigma||^2\geq 4 hp.
\end{equation}
 
The equality sign of \eqref{E:7.b} holds
if and only if the following statements hold:

$(a)$ $N_T^h$ is an open portion of $CP^h(4)$. 

$(b)$ $N_\perp^p$ is a totally geodesic totally real
submanifold.

$(c)$ The immersion is the following composition:
\begin{equation}\label{E:c} N_T^h\times N_\perp^p
\xrightarrow[\text{{\rm
product immersion}}]{\text{{\rm totally geodesic}}}
CP^h(4)\times CP^p(4)\xrightarrow [\text{{\rm
Segre imbedding}}]{\text{$S_{hp}$}}\notag\end{equation} $$
CP^{h+p+hp}(4)\xrightarrow [\text{{\rm
K\"ahlerian}}]{\text{{\rm totally geodesic}}}
 CP^m(4).$$ 
\end{theorem}

\section{$CR$-warped products and partial
Segre $CR$-immersions.}

Let $B$ and $F$ be two Riemannian manifolds {\it of
positive dimensions\/} equipped with Riemannian
metrics $g_B$ and
$g_F$, respectively, and let
$f$ be a positive function on $B$.
Consider the product manifold $B\times F$ with
its natural projections $\pi:B\times F\to B$ and
$\eta:B\times F\to F$. The {\it warped product}
$M=B\times_f F$ is the manifold
$B\times F$ equipped with the Riemannian
structure such that 
\begin{equation}\label{E:warped} ||X||^2=||\pi_*(X)||^2+f^2(\pi(x))
||\eta_*(X)||^2\end{equation}

\noindent for any tangent vector $X\in T_xM$.
Thus, we have $g=g_B+f^2 g_F$. The
function $f$ is called the {\it warping function\/}
 of the warped product (cf. \cite{oneill}). 

It was proved
in \cite[I]{c3} that there does not exist a $CR$-submanifold
 in a K\"ahlerian manifold which is locally the warped
product
$N_\perp\times_f N_T$ of a totally real submanifold $N_\perp$
and a holomorphic submanifold $N_T$. It
was also proved in
\cite[1]{c3} that there do exist many $CR$-submanifold in
complex space forms which are the warped product $N_T\times_f
N_\perp$ of  holomorphic submanifolds $N_T$ and totally real
submanifolds $N_\perp$ with non-constant warping functions
$f$.

A $CR$-submanifold of a
Kaehler manifold $\tilde M$ is called in  \cite{c3} a
{\it $CR$-warped product} if it is the warped product
$N_T\times_f N_\perp $ of a holomorphic
submanifold $N_T$ and a totally real submanifold $N_\perp$,
where $f$ denotes the warping function. 

A $CR$-warped product is called {\it a non-trivial $CR$-warped
product} if its warping function is non-constant.

\begin{example} Let
$\hbox{\bf C}^{m}$ be the complex Euclidean
$m$-space  with
 a natural Euclidean complex coordinate system $\{z_1,\ldots,
z_m\}$.  We put
${\bf C}^m_*={\bf C}^m-\{0\}$. Let $(w_0,\ldots,w_q)$
denote a Euclidean coordinate system on the Euclidean
$(q+1)$-space ${\bf E}^{q+1}$.

Suppose $z: N_T\to {\bf C}^m_*\subset {\bf C}^m$ is a
K\"ahlerian immersion of a  K\"ahlerian manifold of
complex dimension $h$ into
${\bf C}^m_*$ and $w:N_\perp\to S^{q}(1)\subset {\bf
E}^{q+1}$ is an isometric immersion of a Riemannian
$p$-manifold into the unit hypersphere
$S^{q}(1)$ of ${\bf E}^{q+1}$ centered at the origin.

For each natural
number $\alpha\leq h$, we define a map:
\begin{align} \label{E:6.2} C_{hp}^\alpha :N_T\times N_\perp \to
\hbox{\bf C}^{m}\times S^q(1)\to \hbox{\bf C}^{m+\alpha q}
\end{align}  by
\begin{align} \label{E:6.3} C_{hp}^\alpha &
(u,v)=\, 
\big(w_0(v)z_1(u),w_1(v)z_1(u),\ldots,w_q(v)z_1(u),\ldots,\\& 
w_0(v)z_\alpha(u),w_1(v)z_\alpha (u)
,\ldots,w_q(v)z_\alpha(u),z_{\alpha+1}(u),\ldots,z_m(u)\big)
\notag\end{align} 
for $u\in N_T$ and $v\in N_\perp$. Then \eqref{E:6.2} induces an
isometric immersion: 
\begin{align}\label{E:6.4} \hat C_{hp}^\alpha :N_T\times_{f}
N_\perp\to
\hbox{\bf C}^{m+\alpha q}
\end{align}  from the warped product
$N_T\times_{f} N_\perp$ with warping function
$f=\sqrt{\sum_{j=1}^\alpha |z_j(u)|^2}$ into
$\hbox{\bf C}^{m+\alpha q}$ as a $CR$-warped product.

We put $$\hbox{\bf C}^{h}_\alpha=\Bigg
\{(z_1,\ldots,z_{h})\in \hbox{\bf
C}^{h}: \sum_{j=1}^\alpha |z_j|^2\ne 0\Bigg\}.$$

When $z:N_T= {\bf C}^h_*\hookrightarrow  {\bf C}^h$  and 
$w: S^{p}(1)\hookrightarrow  {\bf E}^{p+1}$ are the inclusion maps,
the map \eqref{E:6.4} induces a map:
\begin{align} \label{E} S_{hp}^\alpha :{\bf
C}^{h}_\alpha\times_{f} S^p(1)\to {\bf C}^{h+\alpha p}
\end{align} 
defined by 
\begin{align} \label{6.6} S_{hp}^\alpha
(z,w)&\,=
\big(w_0z_1,w_1z_1,\ldots,w_pz_1,\ldots,
\\& w_0z_\alpha,w_1z_\alpha
,\ldots,w_pz_\alpha,z_{\alpha+1},\ldots,z_h\big)\notag
\end{align} 
for $z=(z_1,\ldots,z_h)\in \hbox{\bf C}^{h}$ and
$w=(w_0,\ldots,w_p)\in S^{p}(1)$ with $\sum_{j=0}^p w_j^2=1$.
The  warping function of $\hbox{\bf C}^{h}_\alpha\times_{f}
S^p(1)$ is given by
$$f=\Bigg\{\sum_{j=1}^\alpha |z_j|^2\Bigg\}^{1/2}.$$  

The map $S_{hp}^\alpha$  is an isometric $CR$-immersion
which is a  $CR$-warped product in ${\bf C}^{h+\alpha p}$. We
simply call such a
$CR$-warped product  in ${\bf C}^{h+\alpha p}$  a {\it
standard partial Segre $CR$-product.}
 \end{example}
 
The standard partial Segre
$CR$-immersion $S_{hp}^1$ is characterized by
the following theorem (see \cite[I]{c3}). 

\begin{theorem} Let $\phi:_T\times_{f}
N_\perp\to {\bf C}^m$ be a non-trivial $CR$-warped product in the
complex Euclidean $m$-space {\bf C}$^m$  with
$\dim_{\hbox{\bf C}} N_T=h$ and $\dim_{\hbox{\bf R}}
 N_\perp=p$.
Then we have:
\begin{enumerate}
\item[(a)] The squared norm of the second fundamental form 
 satisfies the
inequality:
\begin{align}\label{6.5}||\sigma||^2\geq 
2p||\nabla (\ln f)||^2.\end{align}

\item[(b)] The
$CR$-warped product  satisfies the equality
$||\sigma||^2= 
2p||\nabla (\ln f)||^2$ if and only if the following
statements holds:

\item[(b.1)] $N_T$ is an open portion of
$\hbox{\bf C}^h_1$.

\item[(b.2)] $N_\perp$ is an open portion of
the unit $p$-sphere $S^p(1)$. 
 
\item[(b.3)] The warping function is given by $\;f=|z_1|$.

\item[(b.4)] Up to rigid motions of {\bf
C}$^m$, $\phi$ is the standard partial Segre
$CR$-immersion $S_{hp}^1$. More precisely, we have

\begin{align} 
\phi(z,w)&=\big(S^1_{hp}(z,w),0\ldots,0\big)\\&=\big( 
z_1w_0,z_1w_1,\,\ldots, z_1w_p,z_2,\cdots, z_h
,0,\,\ldots,0\big),\notag \end{align}
for \begin{align} z=(z_1,&\ldots,z_h)\in
\hbox{\bf C}^h_1,\quad w=(w_0,\ldots,w_p)\in
S^p(1)\subset \hbox{\bf E}^{p+1}.
\end{align}
\end{enumerate}\end{theorem}

When $\alpha$ is greater than one, the standard partial Segre
$CR$-immersion $S^\alpha_{hp}$ is characterized by the
following.

\begin{theorem}\label{T:4.1} Let $\phi:N_T\times_f N_\perp\to
\hbox{\bf C}^m$ be a 
$CR$-warped product in complex Euclidean $m$-space {\bf
C}$^m$. Then we have
\begin{enumerate}

\item[(1)]  The squared norm of the second fundamental form
of
$\phi$ satisfies
\begin{equation}\label{E:4.1}||\sigma||^2\geq
2p\big\{||\nabla(\ln f)||^2+\Delta(\ln f)\big\}.\end{equation}

\item[(2)]  If the $CR$-warped product satisfies the
equality case of
\eqref{E:4.1}, then we have

\item[(2.i)]  $N_T$ is an open portion of {\bf C}$^h_\alpha$.

\item[(2.ii)]  $N_\perp$ is an open portion of $S^p(1)$.

\item[(2.iii)]  There exists a natural number $\alpha\leq h$
and  a  complex coordinate
system $\{z_1,\ldots,z_h\}$ on {\bf C}$^h$ such that the
warping function $f$ is given by
$$f=\Bigg\{\sum_{j=1}^\alpha z_j\bar z_j\Bigg\}^{1/2}.$$ 

\item[(2.iv)]  Up to rigid motions of {\bf C}$^m$,
$\phi$ is the standard partial Segre $CR$-immersion
$S_{hp}^\alpha$; namely, we have
\begin{align}\label{E:4.2} \phi&\,(z,w)
=\big(S_{hp}^\alpha(z,w),0\ldots,0\big)  \\&=\big(w_0
z_1,\ldots,w_p z_1,\ldots,w_0z_\alpha ,\ldots, w_p z_\alpha
,z_{\alpha +1},\ldots,z_h,0,\ldots,0\big)\notag \end{align}
 for
$z=(z_1,\ldots,z_h)\in
\hbox{\bf C}_\alpha^h$ and $\,w=(w_0,\ldots,w_p)\in
S^p(1)\subset\hbox{\bf E}^{p+1}$.
\end{enumerate}\end{theorem}

\section{Real hypersurfaces as partial Segre
embeddings.}

A contact manifold is an odd-dimensional
manifold $M^{2n+1}$ equipped with  a 1-form $\eta$
such that $\eta\wedge(d\eta)^n\not=0$. 
A curve $\gamma=\gamma(t)$ in a
contact manifold is called a {\it Legendre curve} if
$\eta(\beta'(t))=0$ along $\beta$. 

We put 
$$ S^{2n+1}(c)=\Bigg\{(z_1,\ldots,z_{n+1})
\in\hbox{\bf C}^{n+1}\, :\,
\<z,z\>=\frac{1}{c}>0\Bigg\}.$$
Let $\xi$ be a unit normal vector of $S^{2n+1}(c)$ in {\bf
C}$^{n+1}$. Then $S^{2n+1}(c)$ is a contact manifold endowed
with a canonical contact structure given by the dual 1-form of
$J\xi$, where $J$ is the complex structure on  {\bf
C}$^{n+1}$.

 Legendre curves are known to play an important role in
the study of contact manifolds. For instance, a
diffeomorphism of a contact manifold is a contact
transformation if and only if it maps Legendre curves to
Legendre curves. 

There is a  simple
relationship between Legendre curves 
 and a second order differential
equation obtained in \cite{1997.2}.

\begin{lemma} Let $c$ be a 
positive number and
$z=(z_1,z_2): I\rightarrow S^3(c)
\subset \hbox{\bf C}^2$ be
a unit speed curve, where $I$ is either an open interval or a
circle.    If  $z$ satisfies the following differential
equation:
\begin{align} z''(t)-i\lambda\gamma(t) z'(t)+cz(t)=0
\end{align} for some nonzero
real-valued  function $\lambda$ on
$I$,  then $z=z(t)$ is a Legendre curve in
$S^3(c)$.

Conversely, if $z=z(t)$ is a Legendre curve in
$S^3(c)\subset {\bf C}^2$, then it satisfies the differential
equation {\rm (7.1)}  for some
 real-valued function $\lambda$. \end{lemma}

For real hypersurfaces in complex Euclidean spaces, we have
the following classification theorem.

\begin{theorem} {\rm (\cite{2002})} Let $a$ be a
positive number and
$\gamma(t)=(\Gamma_1(t),\Gamma_2(t))$ be
a unit speed Legendre curve $\gamma:I\to
S^3(a^2)\subset\hbox{\bf C}^2$ defined on
an open interval
$I$.  Then the partial Segre immersion:
\begin{align} \label{0.2} &\hbox{\bf
x}(z_1,\ldots,z_{n},t)=\big(a\Gamma_1(t)z_1,
a\Gamma_2(t) z_1,z_2,\ldots,z_n\big),\quad
z_1\ne 0
\end{align}
defines a $CR$-warped product real hypersurface, $\hbox{\bf
C}_1^{n}\times_{a|z_1|} I$, in {\bf C}$^{n+1}$, where
$$\hbox{\bf C}_1^{n}=\big\{(z_1,\ldots,z_n):z_1\ne 0\big\}.$$

Conversely, up to rigid motions, every  real hypersurface,
which is the warped product
$\,N\times_f I$ of a complex hypersurface $N$  and an open
interval $\,I$, in  $\,{\bf C}^{n+1}$ is either a partial
Segre immersion defined by \eqref{0.2} or a
 product real hypersurface: $\hbox{\bf C}^n\times
C\subset \hbox{\bf C}^{n}\times\hbox{\bf C}^1$ 
of {\bf C}$^{n}$, where $C$ is a real curve in {\bf C}. 
\end{theorem}

The study of real hypersurfaces in non-flat
complex space forms has been
an active field over the past three
decades.
Although these ambient spaces might be
regarded as the simplest after the spaces
of constant curvature, they impose
significant restrictions on the geometry
of their real hypersurfaces. For instance, they
do not admit totally umbilical
hypersurfaces and Einstein hypersurfaces in non-flat complex
space forms. 

Recently, B. Y. Chen and S. Maeda
prove in \cite{cm} the following general result for real
hypersurfaces in non-flat complex space forms.

\begin{theorem} {\rm (\cite{cm})} Every real
hypersurface  in a complex
projective space (or in a complex
hyperbolic space) is locally an irreducible Riemannian
manifold. 

In other words, there do not exist real
hypersurfaces in non-flat complex space forms which are the
Riemannian products of two or more Riemannian manifolds of
positive dimension.
\end{theorem}

On contrast, there do exist many real
hypersurfaces in non-flat complex space forms
which are  warped products. For real hypersurfaces in
complex projective spaces, we have the following
classification theorem.

\begin{theorem} {\rm (\cite{2002})} Suppose that $a$ is a
positive number and
$\gamma(t)=(\Gamma_1(t),\Gamma_2(t))$ is a
unit speed Legendre curve $\gamma:I\to
S^3(a^2)\subset\hbox{\bf C}^2$ defined on an
open interval
$I$.  Let $\hbox{\bf x}:S_*^{2n+1}\times
I\to\hbox{\bf C}^{n+2}$ be the map defined by the partial
Segre immersion:
$$\hbox{\bf
x}(z_0,\ldots,z_{n},t)=\big(a\Gamma_1(t)z_0,
a\Gamma_2(t) z_0,z_1,\ldots,z_{n}\big),\;\;\;
\sum_{k=0}^{n} z_k\bar z_k=1.$$ 
Then 
\begin{enumerate}

\item $\,\hbox{\bf x}$ induces
an isometric immersion
$\psi:S_*^{2n+1}\times_{a|z_0|} I\to
S^{2n+3}$.

\item  The image \
$\psi(S_*^{2n+1}\times_{a|z_0|} I)$
in
$S^{2n+3}$ is invariant under the action of
$U(1)$.

\item the projection
$\psi_\pi:\pi(S_*^{2n+1}\times_{a|z_0|} I)
\to CP^{n+1}(4)$ of $\psi$ via $\pi$ is
a warped product hypersurface
$CP^{n}_0\times_{a|z_0|} I$ in
$CP^{n+1}(4)$.
\end{enumerate}

Conversely, if a  real hypersurface in
$CP^{n+1}(4)$ is a warped product
$N\times_f I$ of a
complex hypersurface $N$ of $CP^{n+1}(4)$ 
and an open interval $I$, then, up to rigid
motions, it is locally obtained in
the way described above via a partial Segre immersion.
\end{theorem}

\section{Complex extensors, Lagrangian submanifolds and Segre
embedding. }

When $\alpha=h=1$, the partial Segre $CR$-immersion 
\begin{align}\label{7.1}
S^1_{1p}:\hbox{\bf C}^*\times S^p(1)\to \hbox{\bf
C}^{p+1}\end{align}  defined in Section 6 is given by
\begin{align}\label{7.2}
S^1_{1p}(z,w)=\big(zw_0,zw_1,\,\ldots, zw_p\big)\end{align} 
for $ z\in
\hbox{\bf C}^*={\bf C}-\{0\}$ and $w=(w_0,\ldots,w_p)\in
S^p(1)\subset \hbox{\bf E}^{p+1}$.

In this section, we discuss the notion of complex extensors
introduced in \cite{1997} which are constructed in a way
similar to \eqref{7.2}. 

Complex extensors are defined in \cite{1997} as 
follows:

Let $z=z(s):I\to \hbox{\bf C}^*\subset {\bf C}$ be a unit speed curve
in the punctured complex plane ${\bf C}^*$ defined on an open
interval $I$. Suppose that
\begin{align} x=(x_1,\ldots,x_m):&\, M^{n-1}\to
S_0^{m-1}(1)\subset {\bf E}^m\\&\; u\mapsto
(x_1(u),\ldots,x_m(u))\notag \end{align}   
is an isometric immersion of a Riemannian $(n-1)$-manifold
$M^{n-1}$ into ${\bf E}^n$ whose image is contained in
$S_0^{m-1}(1)$.  

The {\it complex extensor} of
$x:M^{n-1}\to {\bf E}^m$ via the unit speed curve 
$z:I\to {\bf C}$ is defined to be
 the map:
\begin{align}\label{7.3}
\tau:&\, I\times M^{n-1}\to  {\bf
C}^{m}\\ &
(s,u)\mapsto \big(z(s)x_1(u),\ldots,
z(s) x_m(u)\big)\notag \end{align}
for $s\in I$ and $u\in M^{n-1}$. It was proved in \cite{1997} that
$I\times M^{n-1}$ is isometrically immersed by $\tau$ as a
totally real submanifold in  ${\bf C}^m$.

The complex tensor of the unit hypersphere
$S^{n-1}\hookrightarrow {\bf E^n}$ via a unit speed curve in {\bf C}
is a $SO(n)$-invariant Lagrangian submanifold in {\bf C}$^n$.
In this way, we can construct many $S(n)$-invariant
Lagrangian submanifolds in {\bf C}$^n$.

Now, we recall the definition of  Lagrangian $H$-umbilical 
submanifolds introduced in \cite{1997,1997.2}.  
   
\begin{definition} A {\it Lagrangian $H$-umbilical 
submanifold\/} of a K\"ahlerian manifold is a
 non-totally geodesic Lagrangian submanifold whose second
fundamental form takes the following  form:  
\begin{align}\label{7.5} &
\sigma(e_1,e_1)= \lambda Je_1,\quad \sigma(e_2,e_2)=\cdots =
\sigma(e_n,e_n)=\mu Je_1,\\ & \sigma(e_1,e_j)=\mu Je_j,\quad
\sigma(e_j,e_k)=0,\;
\; j\not=k,
\;\;\;\; j,k=2,\ldots,n\notag\end{align} for some suitable
functions $\lambda$ and $\mu$ with respect to some suitable
orthonormal local frame field $\{e_1,\ldots,e_n\}$. 
\end{definition}
 
The condition \eqref{7.5} is equivalent to the
single condition: 
\begin{align} &
\sigma(X,Y)=\alpha\<\right. JX,\overrightarrow
H\left. \>\<\right. JY,\overrightarrow H\left. \>\overrightarrow H\\
&
\quad+\beta\<\right. \overrightarrow H,\overrightarrow H\left. \>
\{\<X,Y\>\overrightarrow H+
\<\right. JX,\overrightarrow H\left. \>JY +\<\right.
JY,\overrightarrow H\left. \>JX\}\notag
\end{align}  for vectors $X,Y$ tangent to $M$, where
$$\alpha={{\lambda-3\mu}\over{\gamma^3}},\quad
\beta={\mu\over{\gamma^3}},\quad
\gamma={{\lambda+(n-1)\mu}\over n}$$ when $\overrightarrow H\not=0$. 

It is easy to see that non-minimal Lagrangian $H$-umbilical 
submanifold satisfies the following two conditions:
\begin{enumerate}

\item[(a)] $J\overrightarrow H$ is an eigenvector of the
shape operator
$A_{\overrightarrow H}$.

\item[(b)] The restriction of $A_{\overrightarrow H}$ to
$(J\overrightarrow H)^\perp$ is proportional to the identity map.
\end{enumerate}

On the other hand, since the second fundamental
form of every Lagrangian submanifold
satisfies (see \cite{chenogiue1})
$$\<\sigma(X,Y),JZ\>=\<\sigma(Y,Z),JX\>=\<\sigma (Z,X),JY\>$$ for
vectors $X,Y,Z$ tangent to $M$, we know that  Lagrangian
$H$-umbilical  submanifolds are indeed the simplest
Lagrangian submanifolds which satisfy both  Conditions (a)
and (b). Hence, we can regard  Lagrangian $H$-umbilical 
submanifolds as the simplest Lagrangian submanifolds,  next to
the totally geodesic ones.
  
\begin{example} ({\it Whitney's sphere\/}).
Let $w: S^n\rightarrow {\bf C}^n$ be the map defined by
$$w(y_0,y_1,\ldots,y_n)={{1+iy_0}\over {1+y_0^2}}(
y_1,\ldots,y_n),\quad y_0^2+y_1^2+\ldots+y_n^2=1.$$
Then $w$ is a (non-isometric) Lagrangian immersion of the unit
$n$-sphere into
${\bf C}^n$ which is called the {\it Whitney $n$-sphere}.  

The Whitney $n$-sphere is a complex extensor of the inclusion
$\iota:S^{n-1}\rightarrow {\bf E}^n$ via the unit speed curve
$z$ which is an arclength reparametrization of the 
curve  $\varphi:I\rightarrow {\bf C}$ given by
$$f(\varphi)={{\sin \varphi+i\sin \varphi\cos
\varphi}\over{1+\cos^2\varphi}}.$$  
 
Whitney's
$n$-sphere is a Lagrangian $H$-umbilical  submanifold satisfies
\eqref{7.5} with
$\lambda=3\mu$. In fact, up to dilations,
Whitney's $n$-sphere is the only  Lagrangian $H$-umbilical 
submanifold in ${\bf C}^n$ satisfying $\lambda=3\mu$ (see 
\cite{bcm,1997,ros}). 
\end{example} 

\begin{example}  ({\it Lagrangian
pseudo-spheres\/}). For a given real number $b>0$, let
$z:{\bf R}\rightarrow
{\bf C}$ be the unit speed curve given by
$$z(s)={{e^{2bsi}+1}\over{2bi}}.$$ 
With respect to
the induced metric, the complex extensor of
$\iota:S^{n-1}\rightarrow {\bf E}^n$ via this unit speed
curve is a Lagrangian isometric immersion of an open portion
of $S^n(b^2)$ into
${\bf C}^n$. This Lagrangian submanifold is known as a {\it
Lagrangian pseudo-sphere} \cite{1997}.

A Lagrangian pseudo-sphere is a Lagrangian $H$-umbilical 
submanifold satisfying \eqref{7.5} with $\lambda=2\mu$ (see
\cite{1997}). 
\end{example}

Lagrangian pseudo-sphere is characterized by the following.

\begin{theorem}  {\rm \cite{1997}} 
Let $L:M\rightarrow {\bf C}^n$ be a Lagrangian isometric
immersion. Then, up to rigid motions of ${\bf C}^n$,
$L$ is a Lagrangian pseudo-sphere if and only if $L$ is a 
Lagrangian $H$-umbilical  immersion satisfying
\begin{align} & \sigma(e_1,e_1)= 2bJe_1,\;\;\;\;
\sigma(e_2,e_2)=\cdots = \sigma(e_n,e_n)=bJe_1, \\ &
\sigma(e_1,e_j)=bJe_j,\quad \sigma(e_j,e_k)=0,\quad\quad j\not=k,
\;\; j,k=2,\ldots,n,\notag\end{align} for some nontrivial function
$b$ with respect to some suitable orthonormal local frame field. 

Moreover, in this case, $b$ is a nonzero
constant.\end{theorem} 

The following theorem classifies Lagrangian $H$-umbilical
submanifold in ${\bf C}^n$ with $n\geq 3$.

\begin{theorem}  {\rm \cite{1997}} 
Let $n\geq 3$ and $L: M\rightarrow {\bf C}^n$  be a Lagrangian
$H$-umbilical isometric immersion. Then we have:
\begin{enumerate}

\item If $M$ is of constant sectional curvature,
then either $M$ is flat or, up to rigid motions of ${\bf C}^n$,
$L$ is a Lagrangian pseudo-sphere.

\item If $M$ contains no open subset of constant
sectional curvature, then,  up to rigid motions of ${\bf C}^n$, $L$
is a complex extensor of the unit hypersphere of ${\bf E}^n$
via a unit speed curve in ${\bf C}^*$.
\end{enumerate}\end{theorem} 

\begin{remark} Flat Lagrangian $H$-umbilical submanifolds in {\bf
C}$^n$ are not necessary complex extensors (see \cite{1997}).
For the explicit representation formula of flat Lagrangian
$H$-umbilical submanifolds in {\bf C}$^n$, see
\cite{1999}. \end{remark}

\begin{remark} Complex extensors in an indefinite complex
Euclidean space ${\bf C}^n_s$ are introduced and
are investigated in \cite{c20}.   For the relationship between
complex extensors and Lagrangian submanifolds in indefinite
complex Euclidean spaces and their applications, see
\cite{c20}.\end{remark}

\section{Partial Segre $CR$-immersions in complex projective
space. }

 Let
$\hbox{\bf C}^*=\hbox{\bf C}-\{0\}$ and
$\hbox{\bf C}_*^{m+1}=\hbox{\bf
C}^{m+1}-\{0\}$.  Consider the action of
$\hbox{\bf C}^*$ on
$\hbox{\bf C}_*^{m+1}$ defined by
$$\lambda\cdot (z_0,\ldots,z_m)=(\lambda
z_0,\ldots,\lambda z_m)$$ for $\lambda\in {\bf C}^*$, where
$\{z_0,\ldots,z_h\}$ is 
 a natural complex Euclidean coordinate system on
$\hbox{\bf C}^{m+1}_*$.  Let
$\pi(z)$ denote the equivalent class contains
$z$. Then we have a projection: $\pi:
\hbox{\bf C}_*^{m+1}\to \hbox{\bf C}^{m+1}_*/\sim.$ 
It is known that the  set of equivalent classes under $\pi$ is
the complex projective $m$-space $CP^m(4)$.
The coordinate system
$\{z_0,\ldots,z_m\}$ are the homogeneous coordinate system on
$CP^{m}(4)$. Thus,
we have the projection:
$$\pi:\hbox{\bf C}_*^{m+1}\to CP^m(4).$$

For each integer $\alpha$ with $0\leq \alpha\leq h$. We put
$${\bf C}^{h+1}_\alpha=\Bigg\{(z_0,\ldots,z_h)\in
{\bf C}^{h+1}:\sum_{j=0}^\alpha |z_j|^2\ne 0\Bigg\}.$$

Consider
the map:
\begin{align} S_{hp}^\alpha :\hbox{\bf C}_\alpha^{h+1}\times
S^p(1)\to \hbox{\bf C}_*^{h+p+\alpha p+1} \end{align} 
defined by 
\begin{align} \label{5.3} &S_{hp}^\alpha
(z,w)\, =
\big(w_0z_0,w_1z_1,\ldots,\\ w_pz_0,&\ldots,
w_0z_\alpha,w_1z_\alpha
,\ldots,w_pz_\alpha,z_{\alpha+1},\ldots,z_h\big)\notag
\end{align} 
for $z=(z_1,\ldots,z_h)\in \hbox{\bf C}_\alpha^{h}$ and
$w=(w_0,\ldots,w_p)\in S^{p}(1)$.

Since the image of $ S_{hp}^\alpha$ is invariant under the
action of ${\bf C}^*$, the composition:
\begin{align}\label{5.4}\pi\circ S_{hp}^\alpha :\mathbb
C^{h+1}_\alpha &\times S^p(1)\xrightarrow[\text{Segre
embedding}] {\text{$ S_{hp}^\alpha $}}
\mathbb C^{h+p+\alpha p+1}_*\\&\xrightarrow[\text{projection}]
{\text{$\pi$}} CP^{h+p+\alpha p}(4)\notag\end{align}  
induces an isometric
$CR$-immersion: \begin{align}\breve
S_{hp}^\alpha:CP^h_\alpha\times_f S^p(1)\to CP^{h+p+\alpha
p}(4)\end{align}  of the product manifold
$CP^h_\alpha\times S^p(1)$ into  $CP^{h+p+\alpha p}(4)$, where
$CP^h_\alpha$ is the open subset of
$CP^h(4)$ defined by $$CP^h_\alpha=\Bigg\{(z_0,\ldots,z_h)\in
CP^h(4):\sum_{j=0}^\alpha |z_j|^2\ne 0\Bigg\}.$$ 
The  metric on
$CP^h_\alpha\times S^p(1)$ induced via \eqref{5.4} is a warped
product metric with warping function, say $f$. Clearly,
$CP^h_\alpha$ is a non-compact manifold. 
 
We simply called such a
$CR$-warped product immersion $\breve S_{hp}^\alpha$ in 
$CP^{h+p+\alpha p}(4)$ a {\it standard partial Segre
$CR$-immersion in} $CP^{h+p+\alpha p}(4)$.

The standard partial Segre $CR$-immersion $\breve
S^0_{hp}$ is characterized by the following 
theorem. 

\begin{theorem} {\rm ( \cite[II]{c3})} Let $\phi:N_T\times_{f}
N_\perp\to CP^m(4)$  be a $CR$-warped product, where 
$h=\dim_{\hbox{\bf C}} N_T$ and
$p=\dim_{\hbox{\bf R}} N_\perp$. Then we have:
\begin{enumerate}

\item[(a)] The squared norm of the second fundamental form 
 satisfies the
inequality
\begin{align}\label{6.8}||\sigma||^2\geq 
2p||\nabla (\ln f)||^2.\end{align}

\item[(b)]  The
$CR$-warped product  satisfies the  equality
 case of \eqref{6.8} if and only if the following statements
hold:

\item[(b.1)]   $N_T$ is an open portion of 
complex projective $h$-space $CP^h(4)$.

\item[(b.2)]  $N_\perp$ is an open portion of 
a unit $p$-sphere $S^p$. 

\item[(b.3)]  Up to rigid
motions, $\phi$ is  the
composition $\pi\circ S^0_{hp}$,
where $S^0_{hp}$ is the standard partial Segre
$CR$-immersion, that is,
\begin{align}
\phi&(z,w)=\big(S^0_{hp}(z,w),0,\ldots,0\big)\\&=\big(
z_0w_0,\,\cdots ,z_0 w_p,z_1,\ldots,z_h,0,\ldots,0\big),
\notag \end{align} for $z=(z_0,z_1,\ldots,z_h)
\in\hbox{\bf C}^{h+1}_0$ and $w=(w_0,\ldots,w_p)\in
S^p\subset \hbox{\bf E}^{p+1}$, 
and $\pi$ is the natural
projection 
$\pi:\hbox{\bf C}^{m+1}_*\to CP^m(4)$.
\end{enumerate}\end{theorem} 

The  standard partial Segre $CR$-immersions $\breve
S^\alpha_{hp}$ with $\alpha>0$ are characterized by the
following  theorem (see \cite{c6}).

\begin{theorem}\label{T:5.1} Let $\phi:N_T\times_{f}
N_\perp\to CP^m(4)$  be a
$CR$-warped product with $h=\dim_{\hbox{\bf C}} N_T$ and
$p=\dim_{\hbox{\bf R}} N_\perp$. Then we have:
\begin{enumerate}

\item The squared norm of the second fundamental form
of
$\phi$ satisfies the inequality:
\begin{equation}\label{E:5.4}||\sigma||^2\geq
2p\big\{||\nabla (\ln f)||^2+\Delta(\ln
f)\}+4hp.\end{equation}

\item The $CR$-warped product  satisfies the equality
case of \eqref{E:5.4} if and only if the following statements
hold:

\item[(2.a)]  $N_T$ is an open portion of 
complex projective $h$-space $CP^h(4)$.

\item[(2.b)]   $N_\perp$ is an open portion of 
 unit $p$-sphere $S^p$.

\item[(2.c)]  There exists a natural number $\alpha\leq
h$ such that, up to rigid motions,  $\phi$ is
given by 
$\pi\circ  S^\alpha_{hp}$, where $S^\alpha_{hp}$ is the standard
partial Segre $CR$-immersion, that is, 
\begin{align}\label{E:5.5} \phi(z,w)&\,=S^\alpha_{hp}(z,w)
\,\\=\big(w_0 z_0,\ldots,w_p z_0,&\ldots,w_0z_\alpha
,\ldots, w_p z_\alpha ,z_{\alpha +1},\ldots,z_h,0\ldots,0\big)
\notag\end{align} for
$z=(z_0,\ldots,z_h)\in
\hbox{\bf C}_\alpha^{h+1}$ and $\,w=(w_0,\ldots,w_p)\in
S^p\subset\hbox{\bf E}^{p+1}$.
\end{enumerate}\end{theorem}

It follows from Example 7.1 that there exist many
 $CR$-warped products $N_T\times_f N_\perp$ with non-constant
warping function in
$CP^{h+p}$ with
$h=\dim_{\hbox{\bf C}} N_T$ and $p=\dim_{\hbox{\bf R}}
N_\perp$.

On contrast, when $N_T$ is compact, the following theorem
shows that the dimension of the ambient space is at least as
 the dimension of the Segre embedding. 

\begin{theorem} {\rm (\cite{c7})} Let 
$N_T\times_f N_\perp$ with  
$h=\dim_{\bf C} N_T$ and
$p=\dim_{\bf R} N_\perp$ be a 
$CR$-warped product in the complex
projective $m$-space
$CP^m(4)$. If $N_T$ is compact, then we have
\begin{align} m\geq h+p+hp.
\end{align}
\end{theorem}

When the dimension of the ambient space $CP^m$ is $m=h+p+hp$
which is the smallest possible,  we have the
following.

\begin{theorem} {\rm (\cite{c7})} Let $\,N_T\times_f
N_\perp\,$ with $h=\dim_{\bf C} N_T$ and
$p=\dim_{\bf R} N_\perp$ be a $CR$-warped product which is
embedded in $\,CP^{h+p+hp}(4)$. If
$N_T$ is compact, then  $N_T$ is holomorphically isometric
to $CP^h(4)$.
\end{theorem} 

\section{Convolution of Riemannian manifolds.}

The  notion of  convolution
of Riemannian manifolds was introduced in \cite{c4,c5}.
This notion extends the notion of warped products in a
natural way. 

\begin{definition} Let  $(N_1,g_1)$ and
$(N_2,g_2)$ be two Riemannian manifolds and let $f$ and $h$
be two positive differentiable functions on $N_1$ and $N_2$,
respectively. Consider the
symmetric tensor field
${}_hg_1 *_f g_2$ of type (0,2) on
$N_1\times N_2$ defined by
\begin{align}\label{6.1}
{}_hg_1 *_f g_2=h^2 g_1+f^2 g_2+2
fh df\otimes dh.\end{align}

The
symmetric tensor field
${}_hg_1 *_f g_2$ is called the 
{\it convolution of $g_1$ and
$g_2$} via $h$ and $f$. The product manifold
$N_1\times N_2$ together with ${}_hg_1 *_f
g_2$, denoted by ${}_hN_1 \bigstar\,_f N_2$,  is
called a {\it convolution manifold}.  

If $\,{}_hg_1 *_f g_2\,$ is a positive-definite
symmetric tensor, it defines  a Riemannian metric on
$N_1\times N_2$. In this case,  ${}_hg_1 *_f g_2$ is
called a   {\it convolution metric \/} and  the 
convolution manifold
${}_hN_1 \bigstar{_f} N_2$ is  called a {\it
convolution Riemannian  manifold.} \end{definition}

When
$f,h$ are irrelevant,  ${}_hN_1
\bigstar\, _f N_2$ and ${}_hg_1 *_f g_2$ are simply denoted
by
$N_1\bigstar\, N_2$ and $g_1*g_2$, respectively.

 The following
result shows that the  notion  of convolution
manifolds arises very naturally. 

\begin{theorem} {\rm (\cite{c4})} Let $x:(N_1,g_1)\to
\hbox{\bf E}^{n}_*\subset \hbox{\bf E}^n$  and 
$y:(N_2,g_2)\to\hbox{\bf E}^m_*\subset \hbox{\bf
E}^m$ be isometric immersions of Riemannian
manifolds
$(N_1,g_1)$ and
$(N_2,g_2)$ into $\hbox{\bf E}^{n}_*$ and
$\hbox{\bf E}^{m}_*$, respectively. Then the map
\begin{align}\psi\,:&\,N_1\times N_2\to \hbox{\bf
E}^n\otimes \hbox{\bf E}^m=\hbox{\bf E}^{n m};\\&
\;(u,v)\mapsto x(u)\otimes y(v),\quad \quad u\in N_1,\;
v\in N_2, \notag\end{align} gives rise to a convolution
manifold
$N_1\bigstar\, N_2$ equipped  with
\begin{align}{}_{\rho_2} g_1*_{\rho_1} g_2=\rho_2^2
g_1+\rho_1^2 g_2 +2\rho_1\rho_2 d\rho_1\otimes
d\rho_2,\end{align}  where  
$$\rho_1=\Bigg\{\sum_{j=1}^n
x_j^2\Bigg\}^{1/2}\quad \rho_2=\Bigg\{\sum_{\alpha=1}^m
y_\alpha^2\Bigg\}^{1/2}$$ 
denote the distance functions of
$x$ and $y$, and
$$x=(x_1,\ldots,x_n),\quad  y=(y_1,\ldots,y_m)$$ are
Euclidean coordinate systems of $\,\hbox{\bf E}^n$
and $\hbox{\bf E}^m$, respectively.
\end{theorem}

\begin{definition}  Let $\psi: ({}_h N_1 \bigstar{_f} N_2,
{}_h g_1 *{_f} g_2)\to (\tilde M,\tilde g)$ be a map from
a convolution manifold into a Riemannian manifold. Then
the map is said to be {\it isometric} if ${}_h g_1 *{_f} g_2$
is induced from $\tilde g$ via $\psi$, that is, we have
\begin{align}\psi^*\tilde g={}_h g_1 *{_f} g_2.\end{align}
\end{definition} 

\begin{example}  Let $x:(N_1,g_1)\to
{\bf E}^{n}_*\subset {\bf E}^n$ be an isometric
immersion. If $y:(N_2,g_2)\to S^{m-1}(1) \subset
{\bf E}^m$ is an isometric immersion such that $y(N_2)$
is contained in the unit hypersphere
$S^{m-1}(1)$  centered at the origin.
Then the convolution $g_1*g_2$ of $g_1$ and $g_2$  is nothing
but the warped product metric: $g=g_1+|x|^2 g_2$.
\end{example}

\begin{definition} A convolution 
${}_h g_1*\,_f g_2$ of two Riemannian metrics
$g_1$ and $g_2$ is said to be {\it degenerate\/} if
$\det({}_h g_1*_f g_2)=0$ holds identically.
\end{definition}

For $X\in T(N_1)$ 
we denote by $|X|_1$ the
length of $X$ with respect to  metric $g_1$ on
$N_1$. Similarly, we denote by $|Z|_2$ for 
$Z\in T(N_2)$ with
respect to metric 
$g_2$ on $N_2$. 

\medskip
\begin{proposition} {\rm (\cite{c4})} Let ${}_h N_1
\bigstar\,_f N_2$ be the convolution of two Riemannian
manifolds $(N_1,g_1)$ and $(N_2,g_2)$ via $h$ and
$f$. Then ${}_h g_1*_f g_2$ is degenerate if and
only if we have:

{\rm (1)} The length $|\hbox{\rm grad\,} f|_1$
of the gradient
of $f$ on $(N_1,g_1)$ is a nonzero constant, say
$c$.

{\rm (2)} The length $|\hbox{\rm grad\,} h|_2$ of the
gradient of $h$ on $(N_2,g_2)$ is the constant given
by $c^{-1}$, that is, the reciprocal
 of $c$.
\end{proposition}

The following result provides a criterion for a
convolution 
${}_h g_1*_f g_2$ of two Riemannian metrics to be
a Riemannian metric.

\begin{theorem} {\rm (\cite{c4})} Let ${}_h N_1
\bigstar\,_f N_2$ be the convolution of Riemannian
manifolds $(N_1,g_1)$ and $(N_2,g_2)$ via $h$ and
$f$. Then  ${}_h g_1*\,_f g_2$ is a
Riemannian metric on ${}_h N_1
\bigstar\,_f N_2$ if and only if we have
\begin{align}\,|\hbox{\rm grad\,} f|_1\cdot |\hbox{\rm
grad\,} h|_2<1.\end{align}\end{theorem}

\section{Convolutions and Euclidean Segre maps.}

Let ${\bf C}_*^n={\bf C}^n-\{0\}$ and ${\bf
E}_*^m={\bf E}^m-\{0\}$. 
Assume that  $(z_1,\ldots,z_n)$ is a
complex Euclidean coordinate system of {\bf C}$^n$
and $(x_1,\ldots,x_m)$ is a Euclidean coordinate system
on {\bf E}$^m$. 
Suppose that  $z:\hbox{\bf C}^h_*\to\hbox{\bf C}^h$ and
$x:\hbox{\bf E}^p_*\to\hbox{\bf E}^p$ are the inclusion
maps.  

Let $\psi$ be the map:
\begin{align}\label{6.1}\psi:{\bf C}^h_*\times {\bf E}^p_*\to
{\bf C}^{hp}\end{align}
 defined by 
\begin{align}\label{6.2}\psi({z,x})=(z_1x_1,
\ldots, z_1 x_p,\ldots, z_h x_1,\ldots, z_h
x_p)\end{align} for
$z=(z_1,\ldots,z_h)\in \hbox{\bf C}^h_*$ and
$x=(x_1,\ldots,x_p)\in \hbox{\bf E}^p_*$. 
The map \eqref{6.1} is called a {\it Euclidean Segre
map.}

If we put $z_j=u_j+iv_j,  i=\sqrt{-1},$ and
$${{\partial}\over{\partial z_j}}= {1\over 2}\Bigg(
{{\partial}\over {\partial
u_j}}-i{{\partial}\over{\partial v_j}}\Bigg)
$$ for $ j=1,\ldots,h,$
then we obtain from \eqref{6.2} that
\begin{align}\label{6.3}
d\psi\Bigg(\sum_{j=1}^h
z_j{{\partial}\over{\partial
z_j}}\Bigg)=d\psi\Bigg(\sum_{\alpha=1}^h
x_\alpha{{\partial}\over{\partial x_\alpha}}\Bigg) .
\end{align}
Notice that the vector fields $\sum_{j=1}^h
z_j{{\partial}/{\partial
z_j}}$ and $\sum_{\alpha=1}^h
x_\alpha{{\partial}/{\partial
x_\alpha}}$ are nothing but the position vector fields of {\bf
C}$_*^h$ and {\bf E}$_*^p$ in {\bf
C}$^h$ and {\bf E}$^p$, respectively. 

Equation \eqref{6.3} implies that the gradient of
$|z|=\sqrt{\sum_{j=1}^h z_j\bar z_j}$ and 
 the gradient of
$|x|=\sqrt{\sum_{\alpha+1}^p x_\alpha^2}$ are mapped to
the same vector field under $\psi$. 

From \eqref{6.2} and \eqref{6.3} it follow that
$d\psi$ has constant rank $2h+p-1$. Hence 
$\psi(
\hbox{\bf C}^h_*\times
\hbox{\bf E}^p_*)$ gives rise to a
$(2h+p-1)$-manifold, denoted by
\begin{align}{\bf C}^h_*\circledast\, {\bf
E}^p_*,\end{align} which equips  a Riemannian metric induced
from the canonical metric on ${\bf C}^h\otimes {\bf
E}^p$ via $\psi$. 

From \eqref{6.2} we can 
verify that
$\hbox{\bf C}^h_*\circledast\,
\hbox{\bf E}^p_*$ is isometric to the
warped product $\hbox{\bf C}^h_*\times S^{p-1}$ equipped
with the warped product metric
$g=g_1+\rho_1 g_0$, where $\rho_1$ is the
length of the position function of
$\hbox{\bf C}^h_*$ and  $g_0$ is the 
metric of the unit hypersphere $S^{p-1}(1)$.  

If we denote the  vector field in \eqref{6.3} by $V$, 
then $V$ is a tangent vector field of $\hbox{\bf
C}^h_*\circledast \, \hbox{\bf E}^p_*$ with length
$|x|\, |z|$. 
The Riemannian metric on  
$\hbox{\bf C}^h_* \circledast \,   \hbox{\bf
E}^p_*$ is induced from the following
convolution:
\begin{align}\label{6.4}{}_hg_1
*_f g_2=\mu^2 g_1+\lambda^2 g_2+2 \lambda\mu
d\lambda \otimes d\mu,\quad \lambda=|z|,\quad
\mu=|x|.\end{align}

\begin{definition} An isometric map:
\begin{align} \phi:(\hbox{\bf
C}^h_*\times{\bf   E}^p_*,{}_{\rho_2}g_1
*_{\rho_1} g_2)\to ({\bf
C}^{m},\tilde g_0)\end{align}
 is call a {\it $CR$-map} if $\phi$
maps each complex slice ${\mathbb   C}^h_*\times \{v\}$ 
of ${\bf C}^h_*\times {\bf E}^p_*$ into
a complex submanifold  of {\bf C}$^m$ and it maps each real
slice  $\{u\}\times {\mathbb   E}^p_*$ of ${\bf C}^h_*\times
{\bf E}^p_*$  into a totally real
submanifold of {\bf C}$^m$. 
\end{definition}

The following two theorems characterize the Euclidean Segre
maps in very simple ways. These two theorems can be regarded
as the Euclidean versions of Theorem
\ref{T:s3} with $s=2$ and Theorem \ref{T:s4}.

\begin{theorem} {\rm (\cite{c5})}
 Let $\phi: ({\bf   C}^h_*\times {\bf 
E}^p_*,{}_{\rho_2}g_1 *_{\rho_1} g_2)\to {\bf  C}^{m}$ be an
isometric 
$CR$-map. Then we have:
\begin{enumerate}

\item $m\geq hp$.

\item If $m=hp$, then, up to rigid motions of
${\bf C}^{m}$, $\;\phi$ is  the Euclidean Segre map,
that is,
\begin{align}\; \phi&(z,x)=\psi_{z,x}=(z_1x_1,\ldots,z_1
x_p,z_2x_1,\ldots,z_2 x_p,\ldots,z_h x_p).\end{align}.
\end{enumerate}\end{theorem}

\begin{theorem} {\rm (\cite{c5})}
 Let $\phi: ({\bf   C}^h_*\times {\bf
E}^p_*,{}_{\rho_2}g_1 *_{\rho_1} g_2)\to {\bf  C}^{m}$
 be an isometric 
$CR$-map. Then we have:
\begin{equation}\label{E:7.4} ||\sigma||^2\geq
{{(2h-1)(p-1)}\over{|x|^2 |z|^2}}.\end{equation}

 The equality sign of \eqref{E:7.4} holds identically
if and only if, up to rigid motions of $\bf C^m$, 
$\phi$ is obtained from the Euclidean Segre map, that is,
$\phi$ is given by
\begin{align}\;
\phi&(z,x)=(\psi_{z,x},0)\\=(z_1x_1,\ldots&,z_1
x_p,z_2x_1,\ldots,z_2 x_p,\ldots,z_h x_p,0,\ldots,0
).\notag\end{align}\end{theorem}

\section{Skew Segre embedding. }

Motivated from the Segre embedding and the Veronese embeddings, S.
Maeda and Y. Shimizu define in \cite{maeda2} real analytic but
not holomorphic embeddings: 
\begin{align}\label{12.1}
f^n_\alpha:CP^n\Big(\tfrac{2}{\alpha}\Big)\to
CP^{\binom{n+\alpha}  {\alpha}-1}(4)
\end{align}
defined by
\begin{align} & (z_0,\ldots,z_n)\notag \\ \mapsto &
\Bigg(z_0^{\alpha}\bar
z_0^{\alpha},\ldots,\sqrt{\frac{\alpha}{\alpha_0 !\cdots\alpha_n
!}} \sqrt{\frac{\alpha}{\beta_0 !\cdots\beta_n
!}}z_0^{\alpha_0}\cdots z_n^{\alpha_n}\bar z_0^{\beta_0}\cdots
\bar z_n^{\beta_n},\ldots, z_n^{\alpha}\bar z_n^{\beta}\Bigg)
\notag\end{align}
where $\sum_{i=0}^n \alpha_i=\sum_{i=0}^n\beta_i=\alpha$, 
and $(z_0,\ldots,z_n)$ is a homogeneous coordinate system
on $CP^n\big(\frac{2}{\alpha}\big)$.

If $\alpha=1$, the embedding \eqref{12.1} reduces to 
the {\it skew-Segre embedding}: \begin{align} f^n_1:& \, CP^n(2)\to
CP^{n(n+2)}(4);\\ &(z_0,\ldots,z_n)\mapsto (z_i \bar
z_j)_{0\leq i,j\leq n}. \notag
\end{align}

T. Maebashi and S. Maeda prove  in \cite{MM} that the
squared mean curvature function of $\,f_1^n\,$ is constant  equal
to $n^{-1}$.

 Maebashi and Maeda also prove
the following.

\begin{theorem} {\rm (\cite{MM})} The  skew-Segre embedding
$\,f_1^n\,$ is equal to the following composition: 
\begin{align}  CP^n(2) 
\xrightarrow{\text{\rm minimal}}
S^{n(n+2)-1}\Big(\frac{n+1}{n}\Big)  
\xrightarrow {\text{\rm totally umbilical}} \notag\end{align}
\begin{align}  
S^{n(n+2)}(1)\xrightarrow {\text{\rm totally geodesic}}
CP^{n(n+2)}(4).
\notag
\end{align} \end{theorem}

The first part, $\,CP^n(2)\xrightarrow{\text{\rm minimal}}
S^{n(n+2)-1}\big(\frac{n+1}{n}\big)\,$, of the decomposition for the
skew-Segre embedding has already been considered by G. Mannoury
(1867--1956) in 1899 (see
\cite{mann}) 

The skew-Segre embedding $f_1^n:CP^n(2)\to CP^{n(n+2)}(4)$ is
a totally real pseudo-umbilical embedding which has parallel
mean curvature vector. Moreover, the squared norm of the
second fundamental form of the skew-Segre embedding 
is constant. 

For $\alpha=2$,  the embedding \eqref{12.1} reduces to 
\begin{align} &f_2^n:\, CP^n(1)\rightarrow CP^{\binom{n+2} 
{2}-1}(4);\\ (z_0,\cdots,z_n)&\,\mapsto\Big(\cdots,z_i^2\overline
z_k^2,\cdots,z_i^2\overline z_k\overline z_l,\cdots,z_
iz^j\overline z_k^2,\cdots,\sqrt 2z_iz_j\overline
z_k\overline z_l,\cdots\Big). \notag\end{align}

S. Maeda and Y. Shimizu gave in \cite{maeda2} an analogous
decomposition for this  embedding.

\begin{theorem} {\rm (\cite{maeda2})} The  embedding
$\,f_2^n:\, CP^n(1)\rightarrow CP^{\binom{n+2}  
{2}-1}(4)\,$ is equal to the following composition: 
\begin{align}  & \quad CP^n(1)  \xrightarrow{\text{\rm minimal}}
S^{m_1-1}(c_1)\times S^{m_2-1}(c_2)   \xrightarrow
[\text{\rm Clifford embedding}] {\text{\rm natural}}
\\& S^{m_1+m_2-1}(c) \xrightarrow {\text{\rm totally umbilical}}
S^{m_1+m_2}(1) \xrightarrow[\text{\rm totally geodesic}] {\text{\rm
totally real}} CP^{\binom{n+2} 
{2}-1}(4),
\notag\end{align}  
where
\begin{align}  & m_1=n(n+2),\;\;
m_2=\frac{1}{4}n(n+1)^2(n+4),\notag\\&
c_1=\frac{(n+1)(n+3)}{4n},
c_2=\frac{(n+2)(n+3)}{n(n+1)},\notag\\ &\hskip.8in 
c=1+\frac{2}{n(n+3)}.
\notag\end{align}  

\end{theorem}

  They observe that $\,f_2^n\,$ is a pseudo-umbilical totally real
embedding and the mean curvature vector of $f_2^n$ is not parallel
in the normal bundle although the mean curvature is constant.

\section{Tensor product immersions and Segre embedding}

The map \eqref{1.2} which defines the Segre embedding can be
regarded as a tensor product map. Here, we recall the
notions of tensor product maps and direct sum maps (see
\cite{1993,Dec} for details).

Let $V$ and $W$ be two vector spaces over the field of real or
complex numbers. Denote by $V\otimes W$ and $V\oplus W$ the
tensor product and the direct sum of
$V$ and $W$, respectively. Let $\<\;\, ,\;
\>_V$ and $\<\;\, ,\; \>_W$ denote the inner products on $V$ and
$W$ respectively. Then 
$V\otimes W$  and $V\oplus W$ are  inner product spaces with
the inner products defined respectively by
\begin{align} &
\<v\otimes w,x\otimes y\>=\<v,x\>_V\cdot\<w,y\>_W ,\quad 
v\otimes w, x\otimes y
\in V\otimes W, \\& \<v\oplus w,x\oplus y\>=\<v,x\>_V+\<w,y\>_W
,\quad  v\oplus w, x\oplus y
\in V\oplus  W\end{align}  

 By applying
these  algebraic notions, we have the notion of
  {\it tensor product maps} and {\it direct sum maps}:
\begin{align}&f_1\otimes f_2 : M \rightarrow V \otimes
W \\&f_1\oplus f_2 : M \rightarrow V \oplus
W \end{align}
 associated with two given  maps  $f_1: M\rightarrow V$ 
and $f_2: M \rightarrow W$ of a  Riemannian manifold $(M,g)$.
These maps are
 defined by
\begin{align}&(f_1\otimes f_2)(u)=f_1(u)\otimes f_2(u) \in  V
\otimes W,\\& (f_1\oplus f_2)(u)=f_1(u)\oplus f_2(u) \in  V
\oplus W, \quad
 u\in M,
.\end{align}

Similarly, if $f: M\rightarrow V$ 
and $f: N \rightarrow W$ are maps from two Riemannian
manifolds $M$ and $N$ into $V$ and $W$, respectively. Then we
have the {\it box-tensor product map} and the
{\it box-direct sum map}: 
\begin{align}&(f\boxtimes
h)(u,v)=f(u)\otimes h(v),\\& 
(f\boxplus
h)(u,v)=f(u)\oplus h(v)\quad u\in M,\;\; v\in N.
\end{align}

The Segre embedding, the partial Segre immersions, complex
extensors, Euclidean Segre maps, convolutions, as well as
skew-Segre immersions given above can all be expressed
in terms of (box) tensor product maps and (box) direct sum 
maps.

\begin{example} Let $\iota_1:{\bf C}^h_* \to {\bf C}^h$ and
$\iota_2:{\bf C}^p_* \to {\bf C}^p$ be the inclusion maps. Then
the map $S_{hp}$ defined by \eqref{1.2} is nothing but the box
tensor product $\iota_1\boxtimes \iota_2$ of $\iota_1$ and
$\iota_2$.
\end{example}

\begin{example} Let $$z=(z_1,\ldots,z_m):N_T\to {\bf
C}^m_*\subset {\bf C}^m$$ be a K\"ahlerian immersion of a 
K\"ahlerian
$h$-manifold into ${\bf C}^m_*$ and
$$w=(w_0,\ldots,w_p):N_\perp\to S^{q}(1)\subset {\bf E}^{q+1}$$
be an isometric immersion from a Riemannian
$p$-manifold into the unit hypersphere
$S^{q}(1)$. Then the partial Segre $CR$-immersion
$C^\alpha_{hp}$ defined by \eqref{E:6.3} is nothing but the map:
\begin{align} C^\alpha_{hp}=(z^\alpha \boxtimes
w)\boxplus z^\alpha_\perp\,:\; &N_T\times N_\perp\longrightarrow
{\bf C}^{m+\alpha q}\\& (u,v)\mapsto (z^\alpha(u)\otimes
w(v))\oplus z^\alpha_\perp (u),
\notag\end{align}
where 
\begin{align} &
z^\alpha=(z_1,\ldots,z_\alpha):N_T\to {\bf
C}^\alpha_*\subset {\bf C}^\alpha \notag\end{align} and
\begin{align}& z^\alpha_\perp=(z_{\alpha+1},\ldots,z_m):N_T\to {\bf
C}^\alpha_*\subset {\bf C}^{m-\alpha},\notag
\end{align}
\end{example}

\begin{example} Let $z=z(s):I\to \hbox{\bf C}^*\subset {\bf C}$ be a unit speed curve
in the punctured complex plane {\bf C}$^*$ defined on an open
interval
$I$ and let
$$x=(x_1,\ldots,x_m):M^{n-1}\to S_0^{m-1}\subset {\bf E}^m; u\mapsto
(x_1(u),\ldots,x_m(u))$$ be an isometric immersion of a
Riemannian $(n-1)$-manifold $M^{n-1}$ into ${\bf E}^n$ whose image is
contained in the unit hypersphere.
Then the {\it complex extensor} $\tau$ of 
$x$ via the unit speed curve 
$z$ is nothing but the box tensor product $z\boxtimes x$.

\end{example}

\begin{example} Let $z=(z_1,\ldots,z_n):{\bf C}^n_*\to {\bf
C}^n$ be the inclusion map of ${\bf C}^n_*$ and let $$\bar
z=(\bar z_1,\ldots,\bar z_n):{\bf C}^n_*\to {\bf C}^n$$ be the
conjugation of $z:{\bf C}^n_*\to {\bf
C}^n$. Then the map \eqref{12.1} which defines the skew-Segre
embedding is nothing but the box tensor product $z\boxtimes
\bar z$.\end{example}

\begin{example} Let $\iota_n:S^n(1)\to {\bf E}^{n+1}$ be the
inclusion map. Then, up to dilations, the tensor product
immersion: \begin{align}\iota_n\otimes\iota_n:S^n(1)\to {\bf
E}^{(n+1)^2}\end{align}
is nothing but the first standard immersion of $S^n(1)$ (see
\cite{1984}).

\end{example}

\section{Conclusion.}

From the previous sections we know that maps and immersions
constructed in ways similar to the Segre embedding
provide us many nice examples for various important classes of
submanifolds. Moreover, we also see from the previous sections
that such examples have many nice properties.

\begin{example} For the inclusion maps $z: {\bf
C}^h_*\hookrightarrow {\bf C}^h$ and
$w: {\bf C}^p_*\hookrightarrow {\bf C}^p$, the Segre map:
\begin{align} \label{14.1}
S_{hp}=z\boxtimes
w=(z_jw_t)_{0\leq j\leq h,0\leq t\leq p}\end{align}
 gives
rise to a {\it K\"ahlerian immersion} of the
product K\"ahlerian manifold
$CP^h(4)\times CP^p(4)$ into
$CP^{h+p+hp}(4)$.
\end{example}

\begin{example} For a given K\"ahlerian immersion $z: N_T\to {\bf
C}^m_*\hookrightarrow  {\bf C}^m$, a given Riemannian
immersion $w:N_\perp\to S^{q}(1)\hookrightarrow {\bf
E}^{q+1}$, and a natural number $\alpha\leq h$ with
$h=\dim_{\bf C} N_T$, the partial Segre map:
\begin{align} \label{E:14.2} C_{hp}^\alpha :N_T\times N_\perp
\to
\hbox{\bf C}^{m}\times S^q(1)\to \hbox{\bf C}^{m+\alpha q}
\end{align} defined by
\begin{align} \label{E:14.3} C_{hp}^\alpha=(z^\alpha \boxtimes
w)\boxplus z^\alpha_\perp
\end{align} 
gives rise to a $CR$-submanifold. Such construction
provide us many nice examples of $CR$-warped products in complex
Euclidean spaces.
\end{example}

\begin{example} Let $z: {\bf
C}^h_*\hookrightarrow  {\bf C}^h$ and 
$w: S^{p}(1)\hookrightarrow  {\bf E}^{p+1}$ be the inclusion
maps. Then the partial Segre map:
\begin{align} \label{E:14.2} C_{hp}^\alpha :N_T\times N_\perp
\to
\hbox{\bf C}^{m}\times S^q(1)\to \hbox{\bf C}^{m+\alpha q}
\end{align} 
gives rise to the standard partial Segre $CR$-immersion in a
complex projective space. Such standard partial Segre
$CR$-immersions satisfy the equality case of the general
inequality: \begin{align} ||\sigma||^2\geq 2p||\nabla(\ln
f)||^2.\end{align} 
\end{example}

\begin{example}  For a unit speed curve $z=z(s):I\to
\hbox{\bf C}^*\hookrightarrow 
 {\bf C}$ and a spherical isometric immersion
$x:M^{n-1}\to S_0^{m-1}\hookrightarrow  {\bf E}^m$, the complex
extensor $\tau$ of $x$ via $z$ is nothing but the box tensor
product $z\boxtimes x$. Such  box tensor product immersions
provide us many nice examples of  {\it Lagrangian
submanifolds} in complex Euclidean spaces.
\end{example}

\begin{example}  Let $a$ be a
positive number and
$\gamma(t)=(\Gamma_1(t),\Gamma_2(t))$ be
a unit speed Legendre curve $\gamma:I\to
S^3(a^2)\hookrightarrow \hbox{\bf C}^2$ defined on
an open interval
$I$.  Then the partial Segre immersion defined by
\begin{align} \label{14.5}(a\gamma\boxtimes z^1)\boxplus
z^1_\perp:
\hbox{\bf C}_*^{n}\times_{a|z_1|} I\to {\bf C}^{n+1}
\end{align}
defines a warped product real hypersurface in {\bf C}$^{n+1}$.

Conversely, up to rigid motions, every  real hypersurface in
{\bf C}$^{n+1}$ which is the warped product
$N\times_f I$ of a complex hypersurface $N$  and an open
interval $I$ is either the partial Segre immersion given by
\eqref{14.5} or the
 product real hypersurface: $\hbox{\bf C}^n\times
C$ in ${\bf C}^{n+1}$ over a curve $C$ in the complex plane.
\end{example}  

\begin{example}  Let 
$$x:(N_1,g_1)\to{\bf E}^{n}_* \hookrightarrow  \hbox{\bf
E}^n,\quad
y:(N_2,g_2)\to\hbox{\bf E}^m_* \hookrightarrow {\bf E}^m$$ 
be two isometric immersions of Riemannian
manifolds
$(N_1,g_1)$ and
$(N_2,g_2)$ into $\hbox{\bf E}^{n}_*$ and
$\hbox{\bf E}^{m}_*$, respectively. Then the Euclidean Segre
map:
\begin{align} x\boxtimes y:&\, N_1\times N_2\to \hbox{\bf
E}^n\otimes \hbox{\bf E}^m;
\\& \;\; (u,v)\mapsto x(u)\otimes y(v),\quad\quad u\in N_1,\;
v\in N_2, \notag\end{align} 
gives rise to the {\it convolution manifold}
$N_1\bigstar\, N_2$ equipped  with the convolution:
\begin{align}{}_{\rho_2} g_1*_{\rho_1} g_2=\rho_2^2
g_1+\rho_1^2 g_2 +2\rho_1\rho_2 d\rho_1\otimes
d\rho_2,\end{align}  where  $\rho_1=|x|$
and $\rho_2=|y|$ are
the distance functions of $x$ and $y$, respectively.
\end{example}

\begin{example}  Let $z:{\bf C}^n_*\hookrightarrow  {\bf
C}^n$ be the inclusion map of ${\bf C}^n_*$ and let $\bar
z$ be the conjugation of $z:{\bf C}^n_*\to {\bf
C}^n$. Then the skew-Segre map $z\boxtimes \bar z$ gives rise to
a {\it totally real isometric immersion} of $CP^n(2)$ in
$CP^{n(n+2)}(4)$.
\end{example}

\begin{example}  Let $\psi:M\to  CP^m(2)$ be a K\"ahlerian
immersion from a K\"ahlerian manifold $M$ into
$CP^m(2)$ and let $\bar \psi$ be the
conjugation of $\psi:M\to  CP^m(2)$. Then the map:
\begin{align}\psi\boxtimes \bar \psi: M\to CP^{m(m+2)}(4)
\end{align} defined by
\begin{align} (\psi\boxtimes \bar \psi)(u)=\big(z_i(u)\bar
z_j(u)\big)_{0\leq i,j\leq m}\end{align}  is
 a {\it totally real}  immersion immersion from $M$ into
$CP^{m(m+2)}(4)$. 

Such box tensor product immersions $\psi\boxtimes \bar \psi$
provides us a way to construct many examples of {\it totally
real submanifolds} in complex projective spaces.
\end{example}

\end{document}